\documentclass[a4paper]{amsproc}
\usepackage{amssymb}
\usepackage{amscd}
\usepackage{pstricks}
\usepackage[dvips]{graphicx}
\usepackage{pst-3dplot}
\usepackage[all,cmtip]{xy}
\usepackage{tikz}
\usepackage{subfigure}
\usepackage{multirow}
\usepackage{tikz-3dplot}

\theoremstyle{plain}
 \newtheorem{thm}{Theorem}[section]

 \newtheorem{prop}{Proposition}[section]
 \newtheorem{lem}{Lemma}[section]
 \newtheorem{cor}{Corollary}[section]
\theoremstyle{definition}
 \newtheorem{exm}{Example}[section]
 \newtheorem{remark}{Remark}[section]
\newtheorem{dfn}{Definition}[section]

\numberwithin{equation}{section}

\DeclareMathOperator{\diag}{diag}

\DeclareMathOperator{\grad}{grad}
\DeclareMathOperator{\adj}{adj}
\DeclareMathOperator{\proj}{proj}
\DeclareMathOperator{\Cay}{Cay}

\title[]{Complexity of Lasso with Normalized Data, Geometry and Combinatorics}

\author{Vladimir Dragovi\'c and Borislav Gaji\'c}

\address{
Department of Mathematical Sciences  \\
The University of Texas at Dallas   \\
Richardson, TX\\
USA\\
Mathematical Institute\\
of the Serbian Academy of Sciences and Arts\\
Belgrade\\
Serbia}
\email{Vladimir.Dragovic@utdallas.edu}

\address{
Mathematical Institute\\
of the Serbian Academy of Sciences and Arts\\
Belgrade\\
Serbia}
\email{gajab@mi.sanu.ac.rs}


\begin{document}
\begin{abstract}

We observe and prove that the complexity of lasso for normalized data is smaller than for nonnormalized ones by relating this question to extremal combinatorics and algebraic graph theory. We employ a geometric approach to the lasso as a study of the tangency of the level sets of the least square objective function with the polyhedral boundary sets $B(t)$ of the parameters  in $\mathbb R^p$ with the $\ell_1$ norm equal to $t$.  We geometrically derive   closed exact formulae for the solution of the lasso  under the full rank assumption.   We  establish  important general properties of the solutions of the lasso, which are known to be represented as a simple polygonal chain in $\mathbb{R}^p$. Starting from $p=2$ and $p=3$, we show a striking difference in the maximal number of $p$-dimensional orthants a polygonal chain of a lasso solution can intersect in the case of normalized data vs. nonnormalized data.   We prove that  for a lasso solution in the normalized case exists only one coordinate axis that contains a segment of the solution. We prove that in the normalized case, the number $h_{p,2}$ is a general upper bound for the number of segments of a lasso solution intersecting a $p$-dimensional orthant, where $h_{p,r}$ is the maximal number of binary words of length $p$ such that every two words match at least at $r$ spots, $r\le p$. We prove, using spectral graph theory, that $h_{p,2}=2^{p-1}-\binom{p}{p/2}/2$ for $p$ even and $h_{p,2}=2^{p-1}-\binom{p-1}{(p-1)/2}$ for $p$ odd. It was known that for general data the sharp estimate for that number is $2^{p-1}$, which we identify with $h_{p,1}$. We also find an upper bound for the total number of segments of a lasso solution with normalized data in dimension $p$, e.g. $\hat N(2\ell)=\frac{3^{2\ell}+5}{2}
-
4\ell\,
{}_3F_2\!\left(
\begin{matrix}
\frac{1}{2},\;1-\ell,\;\frac{1}{2}-\ell\\[1mm]
\frac{3}{2},\;1
\end{matrix}
;4
\right)$,
 that is significantly less than $(3^p+1)/2$, a well-known sharp upper bound for the nonnormalized case.

\end{abstract}

\maketitle

\section{Introduction}
We study complexity of the lasso as important question from applied statistics. We  provide a geometric approach to the lasso as a regression shrinkage method in linear regression, that allows us to reduce the complexity questions to certain fundamental extremal combinatorial problems. We solve these problems by employing algebraic graph theory. From statistical perspective, the main message of the paper is that the complexity of lasso with normalized data is strictly smaller that the complexity of lasso with nonnormalized data. From pure mathematics perspective, the main result is the calculation of the maximal number of binary words of length $p$ such that the Hamming distance between every two words does not exceed $p-2$. The main beauty of the present research is that these seemingly distant questions and areas of science appear here hand in hand.

Let us recall that the lasso method was introduced in the seminal paper of Robert Tibshrani \cite{Tib} in 1996.
The lasso is a regularized estimation process of minimization of the least square objective function under an $\ell_1$ constraint in $\mathbb R^p$, for $p\ge 2$. By adding a constraint, one can improve the prediction accuracy and the interpretation. The nonconstrained solution with respect to the constrained one has a lower bias, but a larger variance. In the constrained regularization, the coefficients not only shrink, but can be even set to zero. This helps with the selection of predictors, in identifying those with the strongest influence.

Among all possible choices of $\ell_q$ norms to be used as a constraint, the $q=1$ case stands out as the smallest giving a convex problem.
It is also the largest one leading to sparse solutions as the norm decreases.  The sparsity is important in the interpretation problem mentioned
above, and plays a fundamental role in sparse statistical learning and modeling and also in compressed sensing.

Since the lasso is a convex problem with a convex constraint, there are many well developed methods for solving it.
All these existing methods for $p>2$ are based on iterative, numerical algorithms, like the coordinate decent.
Usually, the solution is needed not for a fixed value of the constraint, or, equivalently, of the Lagrange multiplier, but over a sequence of values of the constraint. This method is known as pathwise coordinate descent.

The complexity of lasso was studies in e.g. \cite{MY, Tib2}. Similarly to the situation with the simplex algorithm, see \cite{KM}, it was shown in \cite{MY}, that though empirical complexity is linear in the problem size, the worst case complexity for the total number of segments in a lasso path is exactly $(3^p+1)/2$, where $p$ is the number of predictors.

The normalized data are defined by the assumption that the covariance matrix has all diagonal entries equal $1$. It is a custom in statistics to use normalized data when dealing with lasso computations.

Our initial surprising observation is that for the small dimensions $p=2$ and $p=3$, the complexity of the lasso for normalized data is smaller than the complexity for lasso for general data. The  main aim of the present paper is to prove that for any  $p$ the maximal complexity for normalized data is significantly smaller than the complexity for nonnormalized data and to quantify that difference.

In this paper, our geometric approach to the lasso is a study of the tangency of the level sets of the least square objective function with the polyhedral boundary sets $B(t)$ of the parameters  in $\mathbb R^p$ with the $\ell_1$ norm equal to $t$. Here $t$ decreases from the value $\hat t$, which corresponds to the actual, nonconstrained minimizer of the least square objective function, denoted by $\hat\beta$. To illustrate the method and keep the paper self-contained, we geometrically re-derive   closed exact formulae for the solution of the lasso  under the full rank assumption.   We  establish  important general properties of the solutions of the lasso, which are known to be represented as a simple polygonal chain in $\mathbb{R}^p$ with $\hat\beta$  and the origin as the endpoints.  We prove that the intersection of a chain and the interior of a $p$-dimensional orthant coincides with a segment minus its end points, which belongs to a ray having   $\hat\beta$  as its initial point.   Already in $p=2$ and $p=3$ cases we show a striking difference in the maximal number of $p$-dimensional orthants a polygonal chain of a lasso solution can intersect in the case of normalized data, which are $1$ and $2$ respectively vs. nonnormalized data, which are $2$ and $4$ respectively.  We illustrate the results using  especially crafted examples with hypothetical data. We prove that  for a lasso solution in the normalized case exists only one coordinate axis that contains a segment of the solution. We prove that in the normalized case, the number $h_{p,2}$ is a general upper bound for the number of segments of a lasso solution intersecting the interior of a $p$-dimensional orthant, where $h_{p,r}$ is the maximal number of binary words of length $p$ such that every two words match at least at $r$ spots, $r\le p$. The problem of finding $h_{p,r}$, though a reminiscent of a combinatorial conjecture of  Erd\"os, fully solved by Kleitman in \cite{Kl}, is significantly different than the problem of Erd\"os.

We prove that $h_{p,2}=2^{p-1}-\binom{p}{p/2}/2$ for $p$ even and $h_{p,2}=2^{p-1}-\binom{p-1}{(p-1)/2}$ for $p$ odd. Our proof consists of several steps: we identify $h_{p,2}$ with the clique number of certain graphs; we show that these graphs are complements of a class of Cayley graphs; in the even case, we calculate their independence numbers using the Cvetkovi\'c bound \cite{Cv}; in the odd case, since the Cvetkovi\'c bound is not applicable, we use another combinatorial-geometric construction to calculate the independence number in terms of the independence number of an even case.

In contrast to the normalized case, it was known from \cite{MY} that for general, nonnormalized  data, the sharp estimate for the number of segments that intersect the interior of a $p$-dimensional orthant is $2^{p-1}$, which we identify with $h_{p,1}$.  We also find an upper bound $\hat N(p)$
$$
\hat N(2\ell)=\frac{3^{2\ell}+5}{2}
-
4\ell\,
{}_3F_2\!\left(
\begin{matrix}
\frac{1}{2},\;1-\ell,\;\frac{1}{2}-\ell\\[1mm]
\frac{3}{2},\;1
\end{matrix}
;4
\right),
$$
\noindent and with a similar formula for odd $p=2\ell+1$, for the total number of segments of a lasso solution with normalized data in dimension $p$, which is significantly less than $(3^p+1)/2$, a well-known sharp upper bound for the nonnormalized case from \cite{MY} mentioned above.

The paper is organized as follows. In Section \ref{sec:regression}, we review basic notions about classical linear regression. Following \cite{DG2022a, DG2023, DG2023a}, in Section \ref{sec:equi} we formulate a more geometric approach to classical linear regression (see also \cite{Cr}, \cite{FMF}, \cite{Ro}). In Proposition
\ref{prop:directhyper}, we show that for a given dataset, the set of all hyperplanes   having the same hyperplanar moment of inertia (equal $\mu$)  with respect to the given dataset form an ellipsoid.  We refer to such an ellipsoid as an equimomental ellipsoid. When $\mu$ vary, the equimomental ellipsoids form a homothetic family with the same center, corresponding
to $\hat\beta$, the actual, nonconstrained minimizer of the least square objective function, defined by the linear regression.  In Section \ref{sec:geomlasso}, we reformulate the lasso in geometric terms, as a tangency condition between the polyhedral boundary sets $B(t)$ of the parameters with the $\ell_1$ norm equal to $t$ and the equimomental ellipsoids from a homothetic family with the common center, see Theorem \ref{thm:geomlasso}. The first global geometric properties of the lasso solutions were derived in Proposition \ref{prop:ray} and Theorem \ref{thm:rays}.

We provide explicit closed formulas for the lasso solutions in Theorem \ref{thm:sol1p} for the parts of the solutions which belong to the interiors of the orthants, and in Theorem \ref{th:edgep} for the parts of the solutions which belong to coordinate hyperplanes. An important question of uniqueness of the solution of the lasso, see \cite{Tib2, HTW}, is addressed in Theorem \ref{th:uniquep}. In Proposition \ref{th:noreturnp} we recall ``{\it a no return property}" in general case, that, if a solution
of a lasso leaves the interior of one orthant to enter the interior of another orthant, it never comes to the interior of the former again.

In Section \ref{sec:smallp}, we study the geometric solutions of the lasso in small dimensions in a full detail.
The study in the initial dimension ($k=3$ which is the same as $p=2$) is specifically simple due to Lemma \ref{lemma2d},  valid for normalized data, which leads to the closed, exact formulas for the lasso solutions in this case in Theorem \ref{thm:lasop2}. Lemma \ref{lemma2d} does not generalize to higher dimensions,  even under the normalization assumption, thus the case $k=4$ (which is the same as $p=3$) already carries a lot of complexity of the general case. The explicit and closed formulas for this case are given in Theorem \ref{thm:sol1} and Theorem \ref{th:edge}. Theorem \ref{thm:sol1} describes the lasso solutions that belong to the interior of one of the orthants. Theorem \ref{th:edge} provides the solutions for the lasso which belong to coordinate planes.  We show that the solutions of the lasso that belong to coordinate planes, reduce to an induced lasso problem for $p=2$, which was solved before in Theorem \ref{thm:lasop2}.

In Section \ref{sec:smallp}, we also start the study of the important and intriguing question of the number of orthants the interior of which the lasso solutions with normalized may intersect. In Proposition \ref{p2pass} we prove that for $p=2$ the lasso solution  for a normalized data set can intersect the interior of only one quadrant.  In Proposition \ref{prop:2nonnorm}, we provide necessary and sufficient conditions for $p=2$ and nonnormalized data to have lasso solutions which intersect interiors of two quadrants. Example \ref{exam:2nonnorm} provides an explicit nonnormalized case which admits a lasso solution which intersects interiors of two quadrants.

 Using Example \ref{exam:transfer3d}, we prove that it is possible for a lasso solution
for $p=3$,  for normalized data, to belong to the interiors of two octants.
 In Theorem \ref{prop:threefaces}, we prove that for the lasso solutions
for $p=3$,  for normalized data, it is impossible to belong to the interiors of more than two octants.  For $p=3$, in \cite{MY} an example is provided  for which the lasso solutions with nonnormalized intersect interiors of four  octants.

After developing the geometric methodology and building intuition in the cases of small dimensions in Section \ref{sec:smallp}, in  Section \ref{sec:compmax}, we return to the case of general dimension $p$. We translate the question of an upper bound for the number of segments  of a lasso solution with normalized data that intersect the interior of orthants of maximal dimension $p$ to a purely extremal combinatorial problem. We solve the extremal combinatorial problem, after translating it into a graph theory problem of calculating the independence number of a class of Cayley graphs in Section \ref{sec:comb}.

In the concluding Section \ref{sec:comptotal}, we provide an upper bound  $\hat N(p)$ on the total number of segments for a lasso for the normalized data in any dimension $p$ and give a closed formula for it, showing that it is significantly smaller that $(3^p+1)/2$. For $p=3$ we show that $6$ is the sharp upper bound for the total complexity of lasso with normalized data. We introduce a notion of normalized lasso-type code in dimension $p$. Every lasso solution with normalized data produces one such code. We provide an example of such a code for $p=4$ with $\hat N(4)=19$ labels.

\section{Classical regression and equimomental ellipsoids}\label{sec:regression}
\subsection{Basic notions of classical regression}

We recall basic definitions of classical simple linear regression models. It is assumed that the values $(x^{(i)})_{i=1}^N$ of the predictors are known, fixed values, as for example values set up in advance in the experiment. The values $(y^{(i)})_{i=1}^N$ of the responses are observed values of uncorrelated random variables $Y_i$, $i=1, \dots, N$ with the same variance $\sigma^2$.
A linear relation is assumed between the predictors $x^{(i)}$ and the responses  $(y^{(i)})_{i=1}^N$:
$
EY_i=\alpha +\beta x^{(i)}, \quad i=1, \dots, N.
$
This relation can be rewritten ed as
$
Y_i=\alpha +\beta x^{(i)} + \epsilon_i, \quad i=1, \dots, N,
$
where $\epsilon_i$ are called {\it the random errors.} They are uncorrelated random variables with zero expectation and the same variance $\sigma^2$. In such models the regression is of $Y$ on $x$, i.e. in the vertical direction.

Let a system of $N$ points $(x_1^{(i)}, x_2^{(i)}, \dots, x_k^{(i)})_{i=1}^N$ be given. Define \emph{the centroid} $C$, whose coordinates are the mean values of the coordinates $\bar x_j$ and  define \emph{the variances} $
\sigma^2_{x_j}$:
$$
\bar x_j=\frac{1}{N}\sum_{i=1}^Nx_j^{(i)},\ \  \sigma^2_{x_j}=\frac{1}{N-1}\sum_{i=1}^N(x_j^{(i)}-\bar x_j)^2, \ \  j=1, \dots, k.
$$

 We adopt the full rank assumption for a dataset that there is no affine subspace of a smaller dimension that contains all $N$ points.

Due to the  full rank assumption, all $\sigma^2_{x_j}$, for $j=1, \dots, k$ are non-zero.
Then, \emph{the correlations} $r_{jl}$ and \emph{the covariances} $p_{jl}$ are, for $j,l=1, \dots, k, l\ne j$:
$$
r_{jl}=\frac{p_{jl}}{\sigma_{x_j}\sigma_{x_l}},\, p_{jl}=\frac{1}{N-1}\sum_{i=1}^N(x_j^{(i)}-\bar x_j)(x_l^{(i)}-\bar x_l).
$$
\emph{The covariance matrix} $K$ is a $(k\times k)$ matrix with  diagonal elements
$
K_{jj} = \sigma^2_{x_j}, j=1, \dots, k,
$
and  off-diagonal elements
$
K_{jl}=p_{jl}, j,l=1, \dots, k, l\ne j.
$
The covariance matrix is always symmetric positive semidefinite. However, in this case of the full rank assumption, we have more: $K$ is a positive-definite matrix. In particular, under the full rank assumption, matrix $K$ has the inverse $K^{-1}$ and all its eigenvalues are positive. It is customary in statistics to assume  that the origin of the Cartesian coordinate system coincides with the centroid.  A data set is {\it normalized} if $K_{jj}=1$ for all $j=1, \dots, k$.

\subsection{Equimomental ellipsoids in the classical regression}\label{sec:equi}

Let  a system of $N$ points in $\mathbb{R}^k$ with masses $m_1,...,m_N$ be given under the full rank assumption, where $k\geqslant 2$. The
 point $C$ denotes the center of masses.

Consider a hyperplane $\pi$ in the same space $\mathbb{R}^k$. The hyperplanar moment of inertia for the system of points for the hyperplane $\pi$ is, by definition:
$
J_{\pi}=\sum\limits_{i=1}^Nm_id_i^2,
$
 where $d_i$ is the  perpendicular distance form the $i$-th point to the hyperplane.

In the case of hyperplanar  moments of inertia, a generalization of the Huygens-Steiner
theorem can be formulated as follows:
$J_{\pi}=J_{\pi_1}+md^2,$
where it is assumed that the hyperplane $\pi_1$ contains the center of masses, while $\pi$ is a hyperplane parallel to $\pi_1$ at the distance $d$.
Here $m$ is the total mass of the system of points. The hyperplanar operator of inertia at the point $O$ is defined here as a $k$-dimensional symmetric operator as follows:
$
\langle J_O\mathbf{n_1},\mathbf{n_2}\rangle=\sum\limits_{j=1}^{N}m_j\langle\mathbf{r_j},\mathbf{n_1}\rangle\langle\mathbf{r_j},\mathbf{n_2}\rangle,
$
where $\mathbf{r_j}$ is the radius vector of the point $M_j$.
We see that
$
J_\pi=\langle J_O\mathbf{n},\mathbf{n}\rangle,
$
where $\mathbf{n}$ is the unit vector orthogonal to the hyperplane $\pi$ which contains $O$.

Let a direction $w$ be given. We introduce  the hyperplanar moment of inertia in direction $w$ for the given system of points and for a given hyperplane $\pi\subset \mathbb R^k$, which is not parallel to $w$. We define $J^w_{\pi}$
the hyperplanar moment of inertia in direction $w$ as follows:
$J^w_{\pi}=\sum\limits_{j=1}^Nm_j\hat D_j^2,$
where $\hat D_j$ is the distance between the point $M_j$ and the intersection of the hyperplane $\pi$ with the line parallel to $w$ through $M_j$, for $j=1, \dots, N$. Let  $\textbf{n}$ be the unit vector orthogonal to the hyperplane $\pi$, and $O$ a point contained  in $\pi$. Then the hyperplanar moment of inertia $J^w_{\pi}$ can be rewritten in the form
$
J^w_{\pi}=J_{\pi}/\langle\textbf{w}_0,\textbf{n}\rangle^2=\langle J_O\textbf{n},\textbf{n}\rangle/\langle\textbf{w}_0,\textbf{n}\rangle^2.
$

\begin{prop}\label{prop: HSdirplan}[The hyperplanar directional Huygens-Steiner Theorem]\label{prop:dHST}
 Let the hyperplane $\pi_1$ contain the center of masses $C$ and let $\pi$ be a hyperplane parallel to $\pi_1$. Denote by $J^w_{\pi_1}$ and $J^w_{\pi}$ the corresponding
 directional hyperplanar moments of inertia of a given system of points with the total mass $m$ in the direction $w$. Then
 $J^w_{\pi}=J^w_{\pi_1}+m\hat D^2$,
where $\hat D$ is the distance between the points of intersection of a line parallel to $w$ with the parallel hyperplanes $\pi$ and $\pi_1$.
\end{prop}

Thus, we get a characterization of the center of masses, using Proposition \ref{prop:HSdirplan}:

\begin{cor} Given a direction $w$, the system of points and one hyperplane $\pi$ not parallel to $w$.
Among all the hyperplanes parallel to $\pi$, the least directional hyperplanar moment of inertia in direction $w$ is attained for the hyperplane which contains the center of masses of the system of points.
\end{cor}

 We are going to describe all hyperplanes that have the same directional hyperplanar moment of inertia in direction $w$.
We consider a coordinate system  with the origin at the centroid $C$, and with the direction $w$ coinciding with the axis $Cx_{k}$. This models  the classical linear regression, in vertical direction.

The  hyperplanar operator of inertia at the point $C$, denoted by $J_{C}$, is given by the formula
$
[J_{C}]_{ij}=J_{ij}$, for $i, j=1, \dots, k.$
It is convenient to introduce  matrix $K_1$  as the upper-left $(k-1)\times (k-1)$ submatrix of the matrix $J_C$, and the vector $b$ by:
\begin{equation}\label{eq:K1}
K_1=\left( \begin{matrix}
J_{11}&J_{12}&J_{13}&...&J_{1k-1}\\
J_{12}&J_{22}&J_{23}&...&J_{2k-1}\\
...&&&&&\\
J_{1k-1}&J_{2k-1}&J_{3k-1}&...&J_{k-1k-1}\\
\end{matrix}
\right),
\end{equation}
\begin{equation}\label{eq:b}
b=(J_{1k},...,J_{k-1k})^T.
\end{equation}

Consider a hyperplane $\pi$ which is not parallel to the axis $Cx_k$. An equation of $\pi$ can be  given in the form

\begin{equation}\label{planepi}
x_k=\beta_0-\beta_1x_1-\beta_2x_2-...-\beta_{k-1}x_{k-1}.
\end{equation}
Let us denote by $\beta=(\beta_1,...,\beta_{k-1})$.
\begin{prop}\label{prop:directhyper}
Given a real number $\mu$, let a hyperplane \eqref{planepi} have the directional hyperplanar moment of inertia in direction $Cx_k$ equal to $\mu$. Then the coefficients $\beta_0, \beta_1,...,\beta_{k-1}$ of the hyperplane \eqref{planepi} satisfy the equation
\begin{equation}\label{coefeq}
\langle K_1\beta,\beta\rangle-2\langle b,\beta\rangle+m\beta_0^2+J_{kk}=\mu.
\end{equation}
By varying the moment of inertia $\mu$, the equations \eqref{coefeq} represent the family of homothetic ellipsoids in the space of parameters $\beta_0, \beta_1,...,\beta_{k-1}$.
\end{prop}
\begin{proof}
We calculate the directional hyperplanar moment of inertia of the hyperplane \eqref{planepi}:
$
J_{\pi}^{x_k}=\sum_{i=1}^Nm_id_i^2=\langle K_1\beta,\beta\rangle-2\langle b,\beta\rangle+m\beta_0^2+J_{kk}.
$
This gives the first part of the Proposition. We used here that $C$ is the centroid and consequently, $\sum_{i=1}^N m_ix_j^{(i)}=0$ for $j=1,...,k$.

In the space of parameters $\beta_0, \beta_1,...,\beta_{k-1}$, the equation \eqref{coefeq} defines a family of quadrics that depend on $\mu$.
Since $K_1$ is a positive definite matrix, this is a family of ellipsoids.
\end{proof}

\begin{prop}\label{prop:homotetic} All ellipsoids from the family \eqref{coefeq} have the same center at the point $(0,\hat{\beta})$, where $\hat{\beta}=K_1^{-1}b$ are the coordinates of the hyperplane with the smallest directional hyperplanar moment of inertia in direction $Cx_k$. The formula \eqref{coefeq} can be rewritten in the form
\begin{equation}\label{coefeq1}
\langle K_1(\beta-\hat{\beta}),\beta-\hat{\beta}\rangle=\mu+\langle K_1^{-1}b,b\rangle-m\beta_0^2-J_{kk}.
\end{equation}

\end{prop}
\begin{proof} We have
$
\mu=\langle K_1(\beta-K_1^{-1}b),\beta-K_1^{-1}b\rangle-\langle K_1^{-1}b,b\rangle+m\beta_0^2+J_{kk}.
$
By introducing $\hat{\beta}=K_1^{-1}b$, one gets the formula \eqref{coefeq1}. The coordinates of the vector $\hat{\beta}$ give the coefficients of the hyperplane with the least directional hyperplanar moment of inertia in direction $Cx_k$.
\end{proof}

In the case when all masses are equal to $m_i=\frac{1}{N}$, the matrix $K_1$ coincides with the covariance matrix $K$.
In that case consider the least square objective function
\begin{equation}\label{eq:deff}
f(\beta_0, \dots, \beta_{k-1})=\frac12\sum_{i=1}^N \Big(y^{(i)}-\beta_0-\sum_{j=1}^{k-1}\beta_jx^{(i)}_j \Big)^2.
\end{equation}
In the classical linear regression, the least square estimate $\beta_j^{ls}$
is obtained by minimizing the least square objective function $f=f(\beta_0, \beta_1, \dots, \beta_{k-1})$
over the parameters $(\beta_0, \beta_1, \dots, \beta_{k-1})$:
$$
\min_{\beta_0, \dots, \beta_{k-1}}f(\beta_0, \dots, \beta_{k-1}).
$$

One can observe that
$
f(\beta_0, \beta_1, \dots, \beta_{k-1})=J_{\pi}^{x_k},
$
for $\pi=\beta_0+\sum_{j=1}^{k-1}\beta_jx_j$. Thus the ellipsoids \eqref{coefeq} form the level sets of the function $f$.

From Proposition \ref{prop:homotetic}, we get that
\begin{equation}\label{eq:ls}
\beta^{ls}=\hat{\beta}=K_1^{-1}b,
\end{equation}
which is a well-known formula obtained by the least square method (see for example formula (3.6) from \cite{HTF}).
If the origin is at the centroid, then, according to the hyperplanar directional Huygens-Steiner Theorem, see Proposition \ref{prop:dHST}, $\beta^{ls}_0=0$.
(Proposition \ref{prop:homotetic} shows that the intersection of the family \eqref{coefeq} with the hyperplane $\beta_0=0$ coincides with the family of homothetic ellipsoids of residuals \cite{DG2023}.)

One can always choose coordinates $(x_1, \dots, x_{k-1})$ in which $K_1$ is diagonal so that the formula \eqref{coefeq} takes a simpler form and the direction $w$ coincides with $Cx_k$. In such coordinates,  the  hyperplanar operator of inertia at the point $C$ has a simpler form
$$
J_{C}=\left( \begin{matrix}
J_1&0&0&...&0&J_{1k}\\
0&J_2&0&...&0&J_{2k}\\
...&&&&&\\
0&0&0&...&J_{k-1}&J_{k-1k}\\
J_{1k}&J_{2k}&J_{3k}&...&J_{k-1k}&J_{k}\\
\end{matrix}
\right).
$$

Let us fix $J_{\pi}^{x_k}=\mu$. Then, using that
$
\det(J_C)=J_1J_2...J_{k-1}( J_k-J_{1k}/J_1-...-J_{k-1k}/J_{k-1})
$
one gets that the equation of the family of ellipsoids \eqref{coefeq} has the form
$$
J_1\Big(\beta_1-\frac{J_{1k}}{J_1}\Big)^2+...+J_{k-1}\Big(\beta_{k-1}-\frac{J_{k-1k}}{J_{k-1}}\Big)^2+m\beta_0^2=\mu-\frac{\det(J_{C})}{J_1J_2...J_{k-1}}.
$$
In these coordinates, the center of the ellipsoids is
$
\beta^{ls}=(0, J_{1k}/J_1,\dots,J_{k-1k}/J_{k-1})^T.
$

\section{Geometric reformulation of the classical lasso shrinkage}
\label{sec:geomlasso}

Proposition \ref{prop:directhyper} has the following interpretation. Every hyperplane not parallel to $Cx_k$ corresponds to a point in the dual space with the coordinates $(\beta_0, \beta_1,...,\beta_{k-1})$. In this correspondence, the hyperplanes with the same directional hyperplanar moment of inertia in direction $Cx_k$ equal to $\mu$ form an ellipsoid \eqref{coefeq}. By varying $\mu$ one gets the family of homothetic ellipsoids with the same center at the point $(0, \hat{\beta}_1,...,\hat{\beta}_{k-1})$.

We are going to use this result to reformulate the lasso in geometric terms. Then, using geometric arguments, we are going to solve the lasso in closed formulas.

Suppose that the data set of $N$ points in $\mathbb R^k$ is given: $(x^{(i)}_1,...,x^{(i)}_{k-1},y^{(i)})$, for $i=1,...,N$. The variables $x^{(i)}_j$, $j=1,\dots, k-1$, $i=1,\dots, N$ are the predictor variables, while $y^{(i)}$, $i=1,\dots, N$ are the corresponding responses.

Given $t>0$, the lasso estimate $\beta^{lasso}_j(t)$ is obtained by minimizing the function $f=f(\beta_0, \beta_1, \dots, \beta_{k-1})$ \eqref{eq:deff},
over the parameters $(\beta_0, \beta_1, \dots, \beta_{k-1})$:
$$
\min_{\beta_0, \dots, \beta_{k-1}}f(\beta_0, \dots, \beta_{k-1})\quad\textrm{under the condition that}\quad  \sum_{j=1}^{k-1}|\beta_j|\leqslant t.
$$
If the origin is at the centroid, then $\beta^{lasso}(t)_0=0$ (see \cite{Tib}).

If there is no condition $\sum_{j=1}^{k-1}|\beta_j|\leqslant t$, then the lasso estimate coincides with the least square estimate $\beta_j^{ls}$. From formula \eqref{eq:ls} it follows that $\beta_j^{ls}=\hat{\beta}_j$. In a case when $K_1$ has a diagonal form, $K_1=\diag(J_1,\dots, J_{k-1})$, the formula is $\beta_j^{ls}=-\frac{J_{jk}}{J_j}$, for $1\leqslant j\leqslant k-1$.

Since $\beta^{lasso}(t)_0=0$, we intersect the family of ellipsoids \eqref{coefeq} with the hyperplane $\beta_0=0$. We get the family of ellipsoids $E_{\mu}$ in the space of parameters $(\beta_1,...,\beta_{k-1})$:
\begin{equation}\label{coefeq2}
E_{\mu}: \langle K_1(\beta-\hat{\beta}),\beta-\hat{\beta}\rangle=\mu+\langle K_1^{-1}b,b\rangle-J_{kk}.
\end{equation}

\begin{lem}\label{lem:beta00}
All the points that lie on one fixed ellipsoid from the family \eqref{coefeq2} correspond to the hyperplanes with the same moment of inertia equal to $\mu$, and form a level set of the function $f(0, \beta_1, \dots, \beta_{k-1})$. When  $\mu$ increases, the semiaxes of the ellipsoid increase. An ellipsoid with a  smaller $\mu$ is placed inside all the ellipsoids with a greater $\mu$.
\end{lem}

In what follows, for simplicity, let us denote $k-1=p$.

The lasso condition $\sum_{j=1}^{p}|\beta_j|\leqslant t$ defines a solid $S(t)$ in $\mathbb R^p$, with the boundary
$B(t)=\sum_{j=1}^{p}|\beta_j|=t$. The boundary $B(t)$ is a polyhedral hypersurface in $\mathbb R^p$, which homothetically decreases when $t$ decreases. Let us observe that $B(t)$ is a cross polytope. This is a polytope dual to a hypercube (see \cite{Co, Zi}). For a given $t$ denote by $\mu=\mu(t)$ the smallest value of $\mu$ for which the ellipsoid $E_{\mu}$ touches the boundary $B(t)$.
Thus, we derived the following geometric formulation of the lasso estimates.

\begin{thm}\label{thm:geomlasso}[Geometric formulation of the lasso] When $t>0$ decreases, the lasso estimate $\beta^{lasso}_j(t)$, $j=1, \dots, p$
gives the coordinates of the point of contact of the boundary $B(t)$ with the ellipsoid $E_{\mu(t)}$. When the point of contact $\beta^{lasso}_j(t)$, $j=1, \dots, p$ belongs to the interior of a face of  $B(t)$ of a maximal dimension $p-1$ of $B(t)$ then the gradient of the function defining the ellipsoid $E_{\mu(t)}$ is collinear with the normal to the face, and they are of the opposite orientations. If the point of contact $\beta^{lasso}_j(t)$, $j=1, \dots, p$ belongs to the interior of a face of $B(t)$ of a smaller dimension $j$,  $j<p-1$, then $p-1-j$ of its coordinates are equal to zero and the face belongs to the tangent hyperplane of the ellipsoid at the point of contact.
\end{thm}

Let $\hat{t}$ be the value for $t$ when $B(t)$ contains the center of the ellipsoids: $\sum_{j=1}^{p}|\hat\beta_j|=\hat{t}$. For $t$ close enough to $\hat{t}$ with $t<\hat{t}$,
the hyperplanes, containing the faces of $B(t)$ that are tangent to the ellipsoids, are parallel to each other. Thus, they are orthogonal to the same vector.

\begin{prop}\label{prop:ray} \begin{itemize}
\item [(a)] The locus of the points of tangency of the family of ellipsoids \eqref{coefeq2} with the family of parallel hyperplanes, which contain the faces of  a maximal dimension of $B(t)$, is a segment belonging to a straight  line that contains the center of the ellipsoids $\hat\beta$.

\item [(b)] Let $\pi$ be a hyperplane that contains the center  $\hat\beta$ of the ellipsoids \eqref{coefeq2}. If two parallel hyperplanes both belong to one of the two  half-spaces defined by  $\pi$, then their contact points with the ellipsoids from the family \eqref{coefeq2} belong to the sam ray, which has $\hat\beta$ as the initial point.
\end{itemize}
\end{prop}
\begin{proof} (a) Suppose that two parallel hyperplanes are tangent to two ellipsoids from the family \eqref{coefeq2}. Their contact points are denoted as $\beta_{(1)}$ and $\beta_{(2)}$ respectively. The gradients of the functions defining the ellipsoids at $\beta_{(1)}$ and $\beta_{(2)}$ are proportional. Thus,  one gets $K_1(\beta_{(1)}-\hat\beta)=\lambda K_1(\beta_{(2)}-\hat\beta)$. Since $K_1$ is nonsingular, one gets $\beta_{(1)}-\hat\beta=\lambda(\beta_{(2)}-\hat\beta)$. Hence, $\beta_{(1)}$ and $\beta_{(2)}$ belong to a straight line that contains $\hat\beta$. (b) The case when $\lambda$ is positive corresponds to the situation when both hyperplanes belong to the same half-space, defined by $\pi$.

\end{proof}

If the coefficient of proportionality $\lambda$ from Proposition \ref{prop:ray} is positive,  then the solutions $(\beta_1(t),..., \beta_p(t))$ belong to a ray that has the center of the ellipsoids as the initial point.
The first intersection of this ray with a coordinate hyperplane (denote the coordinate hyperplane of the first intersection by $\beta_i=0$) happens for $t=t_I$. There exists  some $t_I'$ such that $t_I'<t_I$ and for all $t$ such that $t_I'<t<t_I$, the contact point of the ellipsoid with $B(t)$ belongs to the  the $(p-2)$-dimensional face $\beta_i=0$ of $B(t)$. There are two possibilities. The first possibility occurs if for $t=t_I'$ the passing over the face of dimension $p-2$, given by $\beta_i=0$ happens. Then the solutions for $t$ smaller than $t_I'$, but close enough to $t_I'$, belong to another orthant, and now again $\beta_i(t)\ne 0$. Then, the solutions will belong to a segment that is a subset of another ray initiated at the center of the ellipsoids. The second possibility occurs when for $t=t_I'$, some other coefficient $\beta$ becomes zero, say $\beta_j=0$, $j\ne i$. More details and the exact formulas will be presented below, in Section \ref{sec:genp} and in particular see Theorem \ref{th:edgep}. Summarizing, we have

\begin{thm}\label{thm:rays} The solutions of the lasso that belong to the interior of a face of $B(t)$ of a maximal dimension $p-1$, for varying $t$, lie on the set of segments that are subsets of the rays with the center of the ellipsoids as the common initial point. All the solutions from the same orthant lie on one segment of one of the rays. The number of rays coincide with the number of different orthants containing solutions, or equivalently, with the number of passings over $(p-2)$-dimensional faces.
\end{thm}

The presented approach will be used further,  among other things, to find a linear in $p$ formula for an upper bound of the number of possible passings over  $(p-2)$-dimensional faces of the polyhedron  $B(t)$ in an arbitrary dimension $p$, see Theorem \ref{th:numbpass}. A global view on lasso solutions as polygonal chains is given in Theorem \ref{thm:raysfinal}. This interpretation is in the accordance with the fact indicated in \cite{Tib} that $\beta$'s decrease linearly: when $t$ is getting smaller for some constant rate $\Delta t$, $\beta_i$ is getting smaller with the constant rate $\Delta\beta_i$.

\subsection{Geometric solution to the lasso for an arbitrary $p$}\label{sec:genp}

We start from $K_1$ and $b$ as defined in \eqref{eq:K1} and \eqref{eq:b} respectively, for $p=k-1$.
As in \eqref{eq:ls}, we have $\hat \beta=(K_1)^{-1}b$ is the center of the equimomental ellipsoids $E_{\mu}$ from \eqref{coefeq2}.
Let $\hat t=||\hat \beta||_1$ and $\delta_i$ be the vector of signs of $\hat \beta$, so that $\hat t=\sum_{j=1}^p\delta_i^j\hat\beta_j$
and   $\delta_i^j\hat\beta_j\ge 0$, for all $j=1, \dots, p$. Then for $t<\hat t$ close enough to $\hat t$, the solution to the lasso is given by
\begin{equation}\label{eq:solbetap}
\beta(s)=\hat{\beta} +sK_1^{-1}\delta_i.
\end{equation}

The solution $\beta=\beta(s)$ is linear in $s$. A linear relation between $s$ and $t$ is given by $||\beta(s)||_1=t$. Denote by $s_{\alpha}$ the value of the parameter $s$ such that $\beta_{\alpha}(s_{\alpha})=0$, for $\alpha =1, \dots, p$.
Using the formula
$
K_1^{-1}=(1/\Delta)(\adj K_1),
$
where $\Delta =\det (K_1)$ and $\adj K_1$ is the adjugate of $K_1$, or in other words the transpose of the cofactor matrix of $K_1$,
we get
$$
s_{\alpha}=\frac{-\Delta\hat\beta_{\alpha}}{\sum_{j=1}^pK_1^{j\alpha}\delta_i^j}, \quad \alpha=1, \dots, p,
$$
where $K_1^{j\alpha}$ is the $(j, \alpha)$-cofactor of $K_1$.
\begin{thm}\label{thm:sol1p} Let $s_I$ be the one of the $s_{\alpha}$ $\alpha=1, \dots, p$ that is negative, and such that among all the negative values of $s_{\alpha}$ has the minimal absolute value. For all $s$ such that $s_I<s<0$, the solution of the lasso is given by formulas \eqref{eq:solbetap}.
\end{thm}

When $s$ is smaller than $s_I$, the solution of the lasso is the tangent point of the ellipsoid with the face of the boundary polyhedral $B(t(s))$ of dimension $p-1$. One should observe that the formulas \eqref{eq:solbetap} are not valid any more along that face of dimension $p-2$. The exact formulas in this  case will be presented below, see Theorem \ref{th:edgep}.

There are two possible scenarios  here to be considered. In the first one, in the process of decreasing $s$, the contact point reaches a face of dimension $p-3$ of the boundary polyhedral $B(t(s))$. In that case, two out of $p$ coefficients  $\beta$ are equal to zero.

The second scenario materializes if for some $s$, the contact point passes from the face of dimension $p-1$  to an interior of the adjacent face of dimension $p-1$ of the boundary polyhedral $B(t(s))$.

In order to pass to an adjacent face of the maximal dimension $p-1$, for some $s=s_I'$ the ellipsoid should be tangent to that
adjacent face. The coefficients $s_I$ and $s_I'$ should be both negative, hence $s_Is_I'>0$. Here:
$$
s_{\alpha}'=\frac{-\Delta\hat\beta_{\alpha}}{\sum_{j=1, j\ne \alpha}^pK_1^{j\alpha}\delta_i^j-K_1^{\alpha\alpha}\delta_i^{\alpha}}, \quad \alpha=1, \dots, p.
$$
From the conditions $s_{\alpha}s_{\alpha}'>0$, one gets the condition of the transition of the contact point from the interior of one $(p-1)$-dimensional face to the interior of the adjacent $(p-1)$-dimensional face, over their common $(p-2)$-dimensional subface $\beta_{\alpha}=0$:
$
(\sum_{j=1, j\ne \alpha}^pK_1^{j\alpha}\delta_i^j)^2>(K_1^{\alpha\alpha}\delta_i^{\alpha})^2$, for $\alpha=1,\dots,p.$

Theorem \ref{thm:sol1p} gives the solution of the lasso problem for $s$ such that $s_I<s<0$. For $s<s_I$, instead along an interior of a $(p-1)$-dimensional face, the point of tangency moves along one of the $(p-2)$-dimensional faces $\beta_{\alpha}=0$, for some $\alpha\in\{1,\dots, p\}$.
Then, the gradient $\grad E_{\mu}$ of the ellipsoid is perpendicular to the intersection of the $(p-2)$-dimensional face $\beta_{\alpha}=0$ of the boundary polyhedron, that is
\begin{equation}\label{eq:cont2p}
\sum_{j=1, j\ne {\alpha}}^p|\beta_j|=t.
\end{equation}
In other words, the projection of the gradient $\grad E_{\mu}$ of the ellipsoid to the hyperplane $\beta_{\alpha}=0$ is collinear with the vector of the normal to \eqref{eq:cont2p} and the coefficient of the proportionality is negative. The projection of the gradient of the ellipsoid $\grad E_{\mu}=K_1\beta-b$ onto the hyperplane $\beta_{\alpha}=0$ is
$
\proj_{\alpha}\grad E_{\mu}=(K_1)_{(\alpha;\alpha)}\beta_{(\alpha)}-b_{(\alpha)}.
$
The  condition that the projection is perpendicular  to \eqref{eq:cont2p} can be written in the form
\begin{equation}\label{eq:solp}
(K_1)_{(\alpha;\alpha)}\beta_{(\alpha)}-b_{(\alpha)}=s(\delta_i)_{(\alpha)},
\end{equation}
where $(K_1)_{(\alpha;\alpha)}$ denotes the $(p-1)$-dimensional submatrix obtained from the $p$-dimensional covariance matrix  $K_1$ by omitting  the $\alpha$-th row and the $\alpha$-th column.

Summarizing, we get:

\begin{thm}\label{th:edgep} For $s<s_I$, when the ellipsoid touches the boundary polyhedral along the $(p-2)$-dimensional face  $\beta_{\alpha}=0$, the solution of the lasso reduces to the solution of the induced $(p-1)$-dimensional lasso problem in the hyperplane $\beta_{\alpha}=0$ \eqref{eq:solp}. The corresponding $(p-1)\times (p-1)$ covariance matrix  $(K_1)_{(\alpha;\alpha)}$ is obtained from the $p$-dimensional covariance matrix  $K_1$ by omitting  the $\alpha$-th row and the $\alpha$-th column. The center of the equimomental $(p-2)$-dimensional ellipsoids
is $(\hat \beta)_{\alpha}=((K_1)_{(\alpha;\alpha)})^{-1}b_{(\alpha)}$.
\end{thm}

There is an important question of uniqueness of the solution of the lasso, see \cite{Tib2, HTW}.

\begin{thm}\label{th:uniquep}
The solution of the lasso is unique if and only if the full rank condition is satisfied. If the full rank condition is satisfied for the original
$(p\times p)$ covariance matrix  $K_1$, then the full rank condition is going to be satisfied for all  $(p-j)\times (p-j)$ covariance matrices $(K_1)_{(\alpha_1, \alpha_2,\dots, \alpha_j; \alpha_1, \alpha_2,\dots, \alpha_j)}$ of all the induced lasso problems on all the $(p-j-1)$-dimensional faces of the original $(p-1)$-dimensional boundary polyhedron.
\end{thm}

\begin{proof} The proof follows from the fact that a matrix obtained by removing the $j$-th row and the $j$-th column of a positive-definite matrix is positive-definite. As it is well-known, see \cite{Tib2, HTW}, and as it can be easily seen from our analysis above, the full-rank condition is
necessary and sufficient for the uniqueness of the solution to the lasso.

\end{proof}

Now, we are going to derive and summarize some general properties of the lasso solutions.

\begin{prop}\label{th:noreturnp}
If the point of contact passes from the interior of one $(p-1)$-dimensional face of the polyhedron to one of its boundaries, then it can never return to its interior again.
\end{prop}
\begin{proof}
The solutions $\beta_i(s)$ are linear in $s$. If the solution again appears on the face in the same orthant, then the same linear equation in $s$ for the solution is valid: in the first appearance $s$ is in the interval $(a,b)$, and in the second in a disjoint interval $(c,d)$. But, there is some mid-value $\tilde{s}$ between these intervals, for which the same linear expression does not give a solution. This is a contradiction.
\end{proof}

Combining the results of Theorem \ref{thm:rays} and Theorem \ref{th:edgep}, we provide the following global, quite transparent and simple picture of the lasso solutions.

\begin{thm}\label{thm:raysfinal} Each lasso solution $(\beta_1(s),...,\beta_p(s))$, for varying $s$, form a simple polygonal chain in $\mathbb{R}^p$ with $\hat\beta$ and the origin as the endpoints. There are no two segments of the polygonal chain that are parallel. Each segment of this polygonal chain belongs to one of the rays of the original lasso problem or one of the rays of the induced lasso problems of a face of a smaller dimension.
\end{thm}

\section{Geometry and complexity of the lasso in small dimensions}\label{sec:smallp}

 In this section we will study lasso solutions for small $p$. We will study in detail cases $p=2$ and $p=3$ and touch upon the case $p=4$. The case $p=2$ is exceptionally simple under the normalization assumption. It serves as a base of induction for the lasso solutions for $p>2$. Using Lemma \ref{lemma2d},  valid for normalized data, we get the closed, exact formulas for the lasso solutions in $p=2$ case in Theorem \ref{thm:lasop2}.
It is important to note that already in $p=2$ we will show a striking difference in the maximal number
of quadrants the interior of which a lasso solution can intersect in normalized cases (which is $1$, see Proposition \ref{p2pass}) vs. nonnormalized cases (which is $2$, see \cite {MY} and also Proposition \ref{prop:2nonnorm} and Example \ref{exam:2nonnorm}).
We will show that for $p=3$ the maximal number
of octants the interior of which a lasso solution can intersect in normalized cases is $2$, (see Example \ref{exam:transfer3d} and Theorem \ref{prop:threefaces}) while in  nonnormalized cases, it is $4$ (see \cite {MY}).

Theorem \ref{thm:sol1} describes the lasso solutions for $p=3$ that belong to the interior of one of the octants. Theorem \ref{th:edge} provides the solutions for the lasso which belong to coordinate planes.  The solutions of the lasso that belong to coordinate planes, reduce to an induced lasso problem for $p=2$, which was solved before in Theorem \ref{thm:lasop2}.

Using Example \ref{exam:transfer3d}, we prove that it is possible for a lasso solution
for $p=3$,  for normalized data, to belong to the interiors of two octants.
 In Theorem \ref{prop:threefaces}, we prove that for the lasso solutions
for $p=3$,  for normalized data, it is impossible to belong to the interiors of more than two octants.

\subsection{Lasso in $p=2$}

For $p=2$,  we suppose in this subsection, that standardized data are given with the normalized covariance matrix of the predictors
\begin{equation}\label{corrmatrix2d}
\left(\begin{matrix} 1&J_{12}\\
J_{12}&1
\end{matrix}
\right)
\end{equation}
and the standard deviations $s_{x_1}, s_{x_2}, s_y$. Denote by $J_{13}$ and $J_{23}$ the correlations of the response $y$ with the predictors $z_1$ and $z_2$ respectively. Proposition \ref{prop:directhyper} for $p=2$ gives:

\begin{prop} All the planes $y=\beta_0+\beta_1z_1+\beta_2z_2$ with a fixed directional moment of inertia in the direction of the $y$-axis equal to $\mu$  and a fixed value of $\beta_0$, are parameterized as points of the ellipse in the $(\beta_1, \beta_2)$ plane:
\begin{equation}\label{eq:ellipsemu1}
\frac{\mu}{s_y^2}=m{\beta}_0^2+{\beta}_1^2+{\beta}_2^2-2{\beta}_1{J}_{13}
-2{\beta}_2{J}_{23}+2{\beta}_1{\beta}_2{J}_{12}.
\end{equation}
All the ellipses for varying $\mu$ share the same center with the coordinates
\begin{equation}
\label{eq:centorig}
\hat{\beta_1}=\frac{{J}_{13}-{J}_{23}{J}_{12}}{(1-{J}_{12}^2)s_y},\qquad
\hat{\beta_2}=\frac{{J}_{23}-{J}_{13}{J}_{12}}{(1-{J}_{12}^2)s_y}.
\end{equation}
\end{prop}

The Proposition follows from \eqref{coefeq2} and Lemma \ref{lem:beta00}. From it, we see that among all parallel planes, the minimal directional moment of inertia in the direction of the $y$-axis, has the one with $\beta_0=0$. This is in an alignment with the Huygens-Steiner Theorem. Thus, from now on, we may assume that $\beta_0=0$.

The following Lemma shows that for normalized covariance matrices, as in \eqref{corrmatrix2d}, it is convenient to use a basis in which the covariance matrix of the predictors has a diagonal form.

\begin{lem}\label{lemma2d} Let in a basis $e=[e_1,e_2]$, an operator $A:\mathbb{R}^2\to\mathbb{R}^2$ has the matrix
\begin{equation}\label{2dmatrica}
[A]_e=\left(\begin{matrix} 1&c\\
c&1
\end{matrix}
\right).
\end{equation}
Then, the basis $[\tilde{e}_1, \tilde{e}_2]$, in which the matrix of $A$ has a diagonal form is obtained from the basis $[e_1, e_2]$ by rotation by $\frac{\pi}{4}$.
\end{lem}
\begin{proof} The eigenvalues of the matrix \eqref{2dmatrica} are $1+c$ and $1-c$ and the corresponding eigenvectors are $\frac{1}{\sqrt{2}}(1,1)^T$ and $\frac{1}{\sqrt{2}}(-1,1)^T$. Thus,
$\diag(1+c,1-c)=S^T[A]_eS$, where
$$
S=\frac{1}{\sqrt{2}}\left(\begin{matrix} 1&-1\\
1&1
\end{matrix}
\right).
$$
To finish the proof, we observe that $S$ is the matrix of the rotation by $\pi/4$.
\end{proof}

\begin{thm}\label{thm:lasop2} Given $\hat{\beta_1}, \hat{\beta_2}$ from \eqref{eq:centorig} and $t>0$, the solutions ${\beta}_1,{\beta}_2$ of the lasso with $|{\beta}_1|+|{\beta}_2|\leqslant t$ are:

\begin{enumerate}
\item For $t\leqslant\min\{|\hat{\beta}_1+\hat{\beta}_2|, |\hat{\beta}_1-\hat{\beta}_2|\}$, one of the $\beta_i$ vanishes, while the other one is equal to $t$ or $-t$. We have following cases:
\begin{enumerate}
\item\label{1a} if $\hat{\beta}_1+\hat{\beta}_2>t$, and $\hat{\beta}_2-\hat{\beta}_1>t$, then the solution is $(\beta_1,\beta_2)=(0,t)$;
\item\label{1b} if $-\hat{\beta}_1-\hat{\beta}_2>t$, and $\hat{\beta}_2-\hat{\beta}_1>t$, then the solution is $(\beta_1,\beta_2)=(-t,0)$;
\item\label{1c} if $\hat{\beta}_1-\hat{\beta}_2>t$, and $-\hat{\beta}_1-\hat{\beta}_2>t$, then  the solution is $(\beta_1,\beta_2)=(0,-t)$;
\item\label{1d} if $\hat{\beta}_1+\hat{\beta}_2>t$, and $\hat{\beta}_1-\hat{\beta}_2>t$, then  the solution is $(\beta_1,\beta_2)=(t,0)$;
\end{enumerate}

\item For $|\hat{\beta}_1-\hat{\beta}_2|<t<|\hat{\beta}_1+\hat{\beta}_2|$ we have two subcases

\begin{enumerate}
\item\label{2a} if $\hat{\beta_1}+\hat{\beta}_2>t$, then  the solution is
$
{\beta}_1=(t+\hat{\beta}_1-\hat{\beta}_2)/2, \,
{\beta}_2=(t-\hat{\beta}_1+\hat{\beta}_2)/2.
$

\item\label{2b}
if $-\hat{\beta_1}-\hat{\beta}_2>t$, then  the solution is
$
{\beta}_1=(-t+\hat{\beta}_1-\hat{\beta}_2)/2, \,
{\beta}_2=(-t-\hat{\beta}_1+\hat{\beta}_2)/2.
$

\end{enumerate}

\item For $|\hat{\beta}_1+\hat{\beta}_2|<t<|\hat{\beta}_1-\hat{\beta}_2|$ there are two subcases

\begin{enumerate}
\item\label{3a} if $\hat{\beta_2}-\hat{\beta}_1>t$, then  the solution is
$
{\beta}_1=(-t+\hat{\beta}_1+\hat{\beta}_2)/2, \,
{\beta}_2=(t+\hat{\beta}_1+\hat{\beta}_2)/2.
$

\item\label{3b}
if $\hat{\beta_1}-\hat{\beta}_2>t$, then  the solution is
$
{\beta}_1=(t+\hat{\beta}_1+\hat{\beta}_2)/2, \,
{\beta}_2=(-t+\hat{\beta}_1+\hat{\beta}_2)/2.
$

\end{enumerate}

\item For $t>\max\{|\hat{\beta}_1+\hat{\beta}_2|, |\hat{\beta}_1-\hat{\beta}_2|\}$,
$
\beta_1=\hat{\beta}_1,\,
\beta_2=\hat{\beta}_2.
$

\end{enumerate}

\end{thm}
 \begin{proof} We use coordinates $(\tilde{z}_1,\tilde{z}_2, y)$, such that $(\tilde{z}_1,\tilde{z}_2)$ are coordinates in which the $2\times2$ covariance matrix of the predictors \eqref{corrmatrix2d} has a diagonal form. In the new basis, the covariance matrix of the predictors and the response has the form
 $$
  \left( \begin{matrix}\tilde{J}_1&0&\tilde{J}_{13}\\
 0&\tilde{J}_2&\tilde{J}_{23}\\
 \tilde{J}_{13}&\tilde{J}_{23}&1
 \end{matrix}
 \right).
$$
From Lemma \ref{lemma2d}, it follows that the new basis is
$
\tilde{e}_1=(1/\sqrt{2})(e_1+e_2),\,\tilde{e}_2=(1/\sqrt{2})(-e_1+e_2),\, \tilde{e}_3=e_3.
$
Since the first two vectors of the new basis are obtained by the rotation by $\frac{\pi}{4}$ of the first two vectors of the old basis, the coordinate transformation and its inverse are given by
\begin{equation}\label{eq:betabeta}
\tilde{\beta}_1=\frac{1}{\sqrt{2}}(\beta_1+\beta_2),\, \tilde{\beta}_2=\frac{1}{\sqrt{2}}(-\beta_1+\beta_2);\,
{\beta}_1=\frac{1}{\sqrt{2}}(\tilde{\beta}_1-\tilde{\beta}_2),\, {\beta}_2=\frac{1}{\sqrt{2}}(\tilde{\beta}_1+\tilde{\beta}_2).
\end{equation}

In the new coordinates $(\tilde{\beta}_1, \tilde{\beta}_2)$, the region  $|{\beta}_1|+|{\beta}_2|\leqslant t$ obtains a  more convenient form:
\begin{equation}\label{eq:region}
-\frac{\sqrt{2}}{2}t\leqslant \tilde{\beta}_1\leqslant\frac{\sqrt{2}}{2}t, \quad -\frac{\sqrt{2}}{2}t\leqslant\tilde{\beta}_2\leqslant\frac{\sqrt{2}}{2}t.
\end{equation}
In the new coordinates all the planes with the same directional moment of inertia in the $y$-direction equal to $\mu$ (under the assumption that they contain the origin, i.e. $\beta_0=0$) are parameterized as points of the ellipse (see \eqref{eq:ellipsemu1}):
$$
\tilde{J}_1\Big(\tilde{\beta}_1-\frac{\tilde{J}_{13}}{\tilde{J}_1}\Big)^2+\tilde{J}_{2}\Big(\tilde{\beta}_{2}-\frac{\tilde{J}_{23}}{\tilde{J}_{2}}\Big)^2=\frac{\mu}{s_y^2}-\frac{\det(J_{C})}{\tilde{J}_1\tilde{J}_2}.
$$
All the ellipses with varying $\mu$  have the same center $(\hat{\tilde\beta}_1, \hat{\tilde\beta}_2)=(\tilde{J}_{13}/\tilde{J}_1, \tilde{J}_{23}/\tilde{J}_{2})$. In the original coordinates, the center is given in \eqref{eq:centorig}. Since the semiaxes of the ellipse grow when $\mu$ grows, the solution of the lasso is the point of tangency of the ellipse and the boundary  $\tilde{\beta}_i=\pm t\frac{\sqrt{2}}{2}$, for $i=1,2$ of the region \eqref{eq:region}. Moreover, the principal axes of the ellipses are parallel to the coordinate axes.

In the case when $t<\min\{\sqrt{2}|\hat{\tilde{\beta}}_1|,\sqrt{2}|\hat{\tilde{\beta}}_2|\}$, the center of the ellipse $(\hat{\tilde{\beta}}_1, \hat{\tilde{\beta}}_2)$ is in one of the regions $\tilde{I}$, $\tilde{II}$, $\tilde{III}$, $\tilde{IV}$ from Figure \ref{pic2b}. Since the principal axes of the ellipse are parallel to the coordinate axes,  the ellipse that gives the  solution of the lasso contains one of the vertices of the region \eqref{eq:region} i.e. one of the points $(\pm t\frac{\sqrt2}{2},\pm t\frac{\sqrt2}{2})$.
In the original coordinates, the corresponding regions are $I, II, II, IV$, presented on Figure \ref{pic2a}. If, for example the center $(\hat{\tilde{\beta}}_1, \hat{\tilde\beta}_2)$ belongs to the region $\tilde{I}$, then the solution of the lasso is $(\tilde{\beta}_1,\tilde{\beta}_2)=(t\frac{\sqrt2}{2},t\frac{\sqrt2}{2})$. In the original coordinates, the center belongs to the region $I$. Thus, its coordinates satisfy $\hat{\beta}_1+\hat{\beta}_2>t$, and $\hat{\beta}_2-\hat{\beta}_1>t$. Using coordinate transformations \eqref{eq:betabeta}, one gets the solution in the original coordinates: $(\beta_1,\beta_2)=(0.t)$. This is case \ref{1a}. Using the same procedure, one gets cases \ref{1b}, \ref{1c}, \ref{1d}.

In the case when $\sqrt{2}|\hat{\tilde{\beta}}_1|<t<\sqrt{2}|\hat{\tilde{\beta}}_2|$, the ellipse that gives the solution of the lasso is tangent to the one of the segments of the boundary of the region \eqref{eq:region} given by $\tilde{\beta}_2=\pm t\frac{\sqrt2}{2}$. Thus, there are two subcases: the center belongs to the region between $\tilde{I}$ and $\tilde{IV}$, or it belongs to the region between $\tilde{II}$ and $\tilde{III}$. If, for example $(\hat{\tilde{\beta}}_1, \hat{\tilde{\beta}}_2)$ belongs to the region between $\tilde{I}$ and $\tilde{IV}$, then the solution of the lasso is $(\tilde{\beta}_1, \tilde{\beta}_2)=(\hat{\tilde{\beta}}_1, t\frac{\sqrt{2}}{2})$. In the original coordinates, the center belongs to the region between $I$ and $IV$. The coordinates of the centre satisfy $|\hat{\beta}_1-\hat{\beta}_2|<t$ and $\hat{\beta}_1+\hat{\beta}_2>t$. In the original coordinates, the solution is given by the formulae ${\beta}_1=\frac{1}{2}(t+\hat{\beta}_1-\hat{\beta}_2)$, ${\beta}_2=\frac{1}{2}(t-\hat{\beta}_1+\hat{\beta}_2)$. This is the case \ref{2a}. In a  similar manner, one gets the case \ref{2b}.

In a similar way, one can consider the case $\sqrt{2}|\hat{\tilde{\beta}}_2|<t<\sqrt{2}|\hat{\tilde{\beta}}_1|$.

Finally, when $t>\max\{\sqrt{2}|\hat{\tilde{\beta}}_1|,\sqrt{2}|\hat{\tilde{\beta}}_2|\}$, the solution is $\tilde{\beta}_1=\hat{\tilde{\beta}}_1, \tilde{\beta}_2=\hat{\tilde{\beta}}_3$.

\begin{figure}[h] \centering
\subfigure[][]{%
\begin{tikzpicture}\label{pic2a}
\draw[->] (-2,0) -- (3,0) node[anchor= north east] {$\beta_1$ };
\draw[->] (0,-2) -- (0,3) node[anchor= east] {$\beta_2$};
\draw (-1,0)--(0,-1)--(1,0)--(0,1)--(-1,0);
\draw[dashed] (1,0)--(2,1);
\draw[dashed] (0,1)--(1,2);
\draw[dashed] (0,1)--(-1,2);
\draw[dashed] (-1,0)--(-2,1);
\draw[dashed] (-1,0)--(-2,-1);
\draw[dashed] (0,-1)--(-1,-2);
\draw[dashed] (0,-1)--(1,-2);
\draw[dashed] (1,0)--(2,-1);
\path (1,0.2) node {$t$};
\path (0.2,1) node {$t$};
\path (0.5,2) node {$I$};
\path (2,0.5) node {$II$};
\path (-2,0.5) node {$IV$};
\path (-0.5,-2) node {$III$};
\end{tikzpicture}
}
\hspace{1cm}
\subfigure[][]{%
\begin{tikzpicture}\label{pic2b}
\draw[->] (-2,0) -- (3,0) node[anchor= north east] {$\tilde{\beta}_1$ };
\draw[->] (0,-2) -- (0,3) node[anchor= east] {$\tilde{\beta}_2$};
\draw (-0.71,-0.71)--(-0.71,0.71)--(0.71,0.71)--(0.71,-0.71)--(-0.71,-0.71);
\draw[dashed] (0.71,0.71)--(0.71,2);
\draw[dashed] (0.71,0.71)--(2,0.71);
\draw[dashed] (0.71,-0.71)--(2,-0.71);
\draw[dashed] (0.71,-0.71)--(0.71,-2);
\draw[dashed] (-0.71,-0.71)--(-2,-0.71);
\draw[dashed] (-0.71,-0.71)--(-0.71,-2);
\draw[dashed] (-0.71,0.71)--(-2,0.71);
\draw[dashed] (-0.71,0.71)--(-0.71,2);
\path (0.75,0.2) node {\footnotesize{$\frac{t}{\sqrt{2}}$}};
\path (0.2,0.73) node {\footnotesize{$\frac{t}{\sqrt{2}}$}};
\path (2,2) node {$\tilde{I}$};
\path (2,-2) node {$\tilde{II}$};
\path (-2,-2) node {$\tilde{III}$};
\path (-2,2) node {$\tilde{IV}$};
\end{tikzpicture}
}
\caption{The lasso  in the coordinates  $(\beta_1,\beta_2)$ and $(\tilde{\beta}_1,\tilde{\beta}_2)$.}

\end{figure}
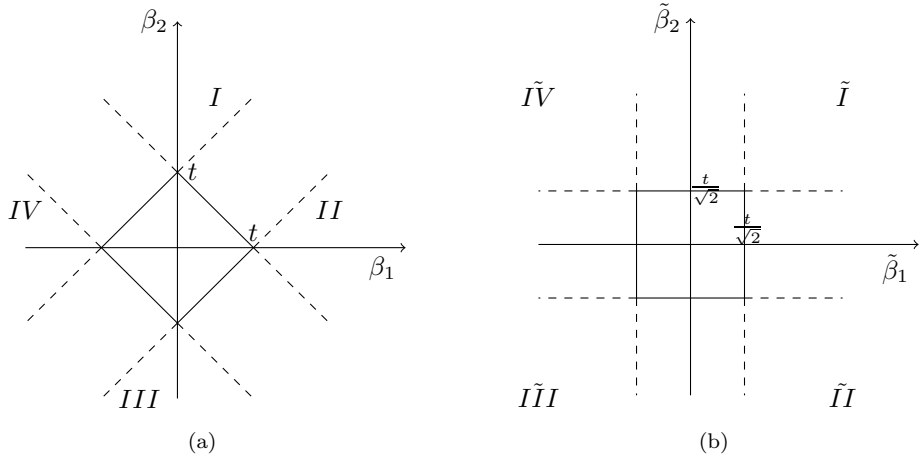

\end{proof}

\begin{figure}[h] \centering

\begin{tikzpicture}\label{pic1p2}
\draw[->] (-2,0) -- (3,0) node[anchor= north east] {$\tilde{\beta}_1$ };
\draw[->] (0,-2) -- (0,3) node[anchor= east] {$\tilde{\beta}_2$};
\draw[] (-0.6,-0.6)--(-0.6,0.6)--(0.6,0.6)--(0.6,-0.6)--(-0.6,-0.6);
\draw[dashed] (-1.3,-1.3)--(-1.3,1.3)--(1.3,1.3)--(1.3,-1.3)--(-1.3,-1.3);
\draw[](2,1) circle [x radius={1.56}, y radius={0.9}];
\draw[dashed](2,1) circle [x radius={0.7}, y radius={0.3}];
\draw[black,fill=gray](2,1) circle [x radius=0.01,y radius=0.01];
\end{tikzpicture}

\caption{The lasso solution for $p=2$ with normalized data.}

\end{figure}

\subsection{ The case $p=3$}

The case $p=2$ is specific  in the normalized setting and, due to Lemma \ref{lemma2d} quite geometrically transparent. There is no analogue of such Lemma for $p\ge 3$  even in the normalized case. Thus, the solution of the lasso in dimension $p=3$ can serve as the baseline for generalizations to get the solutions of the lasso for arbitrary $p$, $p>3$. This is the reason for us to consider the case $p=3$ in great detail before passing to an arbitrary dimension.

The covariance matrix of the predictors is
 \begin{equation}\label{eq:corrmat3}
K_1=\left(\begin{matrix}
1&J_{12}&J_{13}\\
J_{12}&1&J_{23}\\
J_{13}&J_{23}&1
\end{matrix}\right).
\end{equation}
Denote the correlations of the response $y$ with the predictors $z_1$, $z_2$, and $z_3$ by $J_{14}$, $J_{24}$, and $J_{34}$ respectively.

One can suppose again that $\beta_0=0$. All the points of the ellipsoid
\begin{equation}\label{eq:elip3}
\frac{\mu}{s_y^2}={\beta}_1^2+{\beta}_2^2+\beta_3^2-2{\beta}_1{J}_{14}
-2{\beta}_2{J}_{24}-2\beta_3J_{34}+2{\beta}_1{\beta}_2{J}_{12}+2{\beta}_1{\beta}_3{J}_{13}+2{\beta}_2{\beta}_3{J}_{23}
\end{equation}
correspond to the hyperplanes with the moment of inertia equal to $\mu$. The condition
\begin{equation}\label{eq:region3}
|\beta_1|+|\beta_2|+|\beta_3|\leqslant t
\end{equation}
in the space $(\beta_1,\beta_2,\beta_3)$ represents an octahedron, were we suppose that $t$ is fixed.
Since the semiaxes of the ellipsoid grow as $\mu$ grows, the solution of the lasso is the point of tangency of the ellipsoid from the homothetic family \eqref{eq:elip3} and the boundary of the region \eqref{eq:region3} with the corresponding $t=t(\mu)$. All the ellipsoids have the same center $\hat{\beta}$, see \eqref{eq:ls}.
There exists a unique value of $t$, denoted as $\hat{t}$, such that one face of the octahedron \eqref{eq:region3} for $t=\hat {t}$ contains the center of the ellipsoids $\hat{\beta}$.
 Let $\hat{\delta}_i^j$ be $+1$ for $\hat{\beta}_j$ positive, and let it be $-1$ for $\hat{\beta}_j$ negative. Thus
\begin{equation}\label{hatt}
\hat{\delta}_i^1\hat{\beta}_1+\hat{\delta}_i^2\hat{\beta}_2+\hat{\delta}_i^3\hat{\beta}_3=\hat{t}
\end{equation}
is satisfied.
There are $2^3=8$ in total possible choices  of the signs that correspond to $2^3$ faces of the octahedron. For all $t<\hat{t}$ that are close enough to $\hat{t}$, the solution $\beta=\beta(t)$ of the lasso has the same signs as $\hat{\beta}$. The point of contact of the ellipsoid with the boundary of the octahedron can be obtained from the condition that the gradient  to the ellipsoid and the gradient to the corresponding face  of the octahedron are collinear. The coefficient of
proportionality $s$ of two gradients is  negative. One gets the equation for $\beta$:
$
K_1\beta-(J_{14}, J_{24}, J_{34})^T=s(\delta_i^1, \delta_i^2, \delta_i^3)^T.
$

By applying the inverse of matrix $K_1$, for values of $t$ smaller than but close enough to $\hat{t}$ one gets
\begin{equation}\label{eq:solbeta3}
\left(\begin{matrix}\beta_1\\ \beta_2\\ \beta_3\end{matrix}\right)=\left(\begin{matrix}\hat{\beta}_1\\ \hat{\beta}_2\\ \hat{\beta}_3\end{matrix}\right)+
\frac{1}{\Delta}\left(\begin{matrix}
1-J_{23}^2&J_{13}J_{23}-J_{12}&J_{12}J_{23}-J_{13}\\
J_{13}J_{23}-J_{12}&1-J_{13}^2&J_{12}J_{13}-J_{23}\\
J_{12}J_{23}-J_{13}&J_{12}J_{13}-J_{23}&1-J_{12}^2
\end{matrix}\right)\left(\begin{matrix} s\delta_i^1\\ s\delta_i^2\\ s\delta_i^3\end{matrix}\right),
\end{equation}
where
\begin{equation}\label{delta}
\Delta=1+2J_{12}J_{13}J_{23}-J_{12}^2-J_{13}^2-J_{23}^2,
\end{equation}
and $\hat{\beta}$ is given by \eqref{eq:ls}.
The solution is linear in $s$. Additionally, a  linear relation between $s$ and $t$ can be obtained from
$
\hat{\delta}_i^1{\beta}_1+\hat{\delta}_i^2{\beta}_2+\hat{\delta}_i^3{\beta}_3=t.
$

Let us scalar multiply the relation \eqref{eq:solbeta3} by the vector $\delta_i=(\hat{\delta}_i^1,\hat{\delta}_i^2,\hat{\delta}_i^3)$. One gets the following relation between $s$ and $t$;
\begin{equation}\label{tprekos}
t=\hat{t}+s\langle K_1^{-1}\delta_i,\delta_i\rangle,
\end{equation}
where $\hat{t}$ is given by \eqref{hatt}. Since $t<\hat{t}$, and since $K_1$ and its inverse are both positive definite matrices, one gets
\begin{prop} For any given $t$ such that $0<t<\hat{t}$, the value of $s$ that corresponds to the lasso solution that satisfies \eqref{eq:region3} is negative: $s<0$.
\end{prop}

We now deal with the question to find the values of $s$, or equivalently of $t$, for which the formula \eqref{eq:solbeta3} is valid.
Let us denote by $s_i$ the solution of the linear in $s$ equation $\beta_i=0$, for $i=1, 2, 3$. Then
\begin{equation}\label{eq:sfor3}
\begin{aligned}
s_1&=\frac{-J_{14}(1-J_{23}^2)-(J_{13}J_{23}-J_{12})J_{24}-(J_{12}J_{23}-J_{13})J_{34}}{(1-J_{23}^2)\hat{\delta}_i^1+(J_{13}J_{23}-J_{12})\hat{\delta}_i^2+(J_{12}J_{23}-J_{13})\hat{\delta}_i^3},\\
s_2&=\frac{-J_{24}(1-J_{13}^2)-(J_{13}J_{23}-J_{12})J_{14}-(J_{12}J_{13}-J_{23})J_{34}}{(1-J_{13}^2)\hat{\delta}_i^2+(J_{13}J_{23}-J_{12})\hat{\delta}_i^1+(J_{12}J_{13}-J_{23})\hat{\delta}_i^3},\\
s_3&=\frac{-J_{34}(1-J_{12}^2)-(J_{13}J_{12}-J_{23})J_{24}-(J_{12}J_{23}-J_{13})J_{14}}{(1-J_{12}^2)\hat{\delta}_i^3+(J_{12}J_{13}-J_{23})\hat{\delta}_i^2+(J_{12}J_{23}-J_{13})\hat{\delta}_i^1}.
\end{aligned}
\end{equation}
\begin{thm}\label{thm:sol1} Let $s_I$ be the one of the $s_1, s_2, s_3$  that is negative, and such that among all the negative values of $s_1, s_2, s_3$ has the minimal absolute value. For all $s$ such that $s_I<s<0$, the solution of the lasso is given by formulas \eqref{eq:solbeta3}.
\end{thm}
\begin{proof}
The case $s=0$ gives the solution $\beta=\hat{\beta}$. The formulas \eqref{eq:solbeta3} are valid until the point of tangency of the ellipsoid with the face of the octahedron does not reach one edge of
the octahedron. When the point of tangency touches an edge of the octahedron, one of the $\beta_i$ becomes equal to zero. From \eqref{eq:solbeta3} one gets three possible values of $s$ given by
\eqref{eq:sfor3}, each corresponding to the moment of reaching one of the three edges of the face. The point of tangency, actually touches the edge corresponding to $s$ with the minimal absolute value among the values $s_1, s_2, s_3$  that are negative. This finishes the proof.
\end{proof}

When $s$ is smaller than $s_I$, the solution of the lasso is the tangent point of the ellipsoid with the edge. One should observe that the formulas \eqref{eq:solbeta3} are not valid any more in this case. The exact formulas in this  case of tangency along an edge will be presented below, see Theorem \ref{th:edge}.

In dimension $p=2$, we proved that the principal axes of the ellipses are parallel to the edges of the boundary  squares $|\beta_1|+|\beta_2|=t$. This fact  shows that the point of tangency cannot pass from an interior of one edge of the boundary  to the interior of another edge  in the normalized case.  The statement and a proof of this fact will be given in Proposition \ref{p2pass}.

In dimension $p=3$ we show that it is possible   in the normalized case for the point of tangency to pass from the interior of one face to the interior of another face, see Example  \ref{exam:transfer3d}. We will give   necessary and sufficient conditions for $p=3$  in the normalized case when it is possible for the point of tangency to pass from the interior of one face to the interior of another face in Proposition \ref{prop:transfer3d}.

Let us suppose that $s_I=s_1$, or that for $s=s_1$ the solution of the lasso is on the edge $\beta_1=0$. In order to pass to the adjacent face, for some $s=s_1'$ the ellipsoid should be tangent to that
adjacent face. The coefficients $s_1$ and $s_1'$ should be both negative, hence $s_1s_1'>0$. The first component of the normal to the adjacent face changes the sign with respect to the normal of the previous face. Hence, the ellipsoid is tangent to the adjacent face, if
$$
\left(\begin{matrix}\beta_1\\ \beta_2\\ \beta_3\end{matrix}\right)=
\frac{1}{\Delta}\left(\begin{matrix}
1-J_{23}^2&J_{13}J_{23}-J_{12}&J_{12}J_{23}-J_{13}\\
J_{13}J_{23}-J_{12}&1-J_{13}^2&J_{12}J_{13}-J_{23}\\
J_{12}J_{23}-J_{13}&J_{12}J_{13}-J_{23}&1-J_{12}^2
\end{matrix}\right)\left(\begin{matrix} J_{14}-s\delta_i^1\\ J_{24}+s\delta_i^2\\ J_{34}+s\delta_i^3\end{matrix}\right).
$$
The value $s_1'$ is obtained in the limit case, when $\beta_1=0$. One gets:
\begin{equation}\label{s1prim}
s_1'=\frac{-J_{14}(1-J_{23}^2)-(J_{13}J_{23}-J_{12})J_{24}-(J_{12}J_{23}-J_{13})J_{34}}{-(1-J_{23}^2)\hat{\delta}_i^1+(J_{13}J_{23}-J_{12})\hat{\delta}_i^2+(J_{12}J_{23}-J_{13})\hat{\delta}_i^3}\\
\end{equation}
From the conditions $s_1s_1'>0$, one gets:
$
[(J_{13}J_{23}-J_{12})\hat{\delta}_i^2+(J_{12}J_{23}-J_{13})\hat{\delta}_i^3]^2-(1-J_{23})^2>0,
$
or
$
(1-\hat{\delta}_i^2\hat{\delta}_i^3J_{23})^2[(\hat{\delta}_i^2J_{12}+\hat{\delta}_i^3J_{13})^2-(1-\hat{\delta}_i^2\hat{\delta}_i^3J_{23})^2]>0.
$
Thus, we have
\begin{equation}\label{eq:cond1}
(\hat{\delta}_i^2J_{12}+\hat{\delta}_i^3J_{13})^2>(1+\hat{\delta}_i^2\hat{\delta}_i^3J_{23})^2.
\end{equation}
In the same manner, one gets a similar condition for passing from one face to the adjacent along the edge $\beta_2=0$:
\begin{equation}\label{eq:cond2}
(\hat{\delta}_i^1J_{12}+\hat{\delta}_i^3J_{23})^2>(1+\hat{\delta}_i^1\hat{\delta}_i^3J_{13})^2,
\end{equation}
 and a condition for passing from one face to the adjacent along the edge $\beta_3=0$:
 \begin{equation}\label{eq:cond3}
 (\hat{\delta}_i^1J_{13}+\hat{\delta}_i^2J_{23})^2>(1+\hat{\delta}_i^1\hat{\delta}_i^2J_{12})^2.
\end{equation}
 Each of the equations \eqref{eq:cond1}, \eqref{eq:cond2}, \eqref{eq:cond3}, gives actually two conditions by a choice of $\hat{\delta}_i$. The two conditions correspond to passing along consecutive edges of the quadrilateral
 $\beta_{j}=0$, for each $j=1,2,3$. Summarizing, we have the following theorem.

\begin{thm}\label{th:transitions} If the point of tangency of the ellipsoid and the octahedron that gives the lasso solution can pass from  the interior of one face to the  interior of an adjacent face  over the edge $\beta_j=0$ for some $j=1,2,3$, then the corresponding of the following conditions are satisfied:
\begin{equation}\label{eq:usloviprelaz}
\begin{aligned}
&(J_{12}+J_{13})^2-(1+J_{23})^2>0, \, (J_{12}-J_{13})^2-(1-J_{23})^2>0,\quad \textrm{for $\beta_1=0$}, \\
&(J_{12}+J_{23})^2-(1+J_{13})^2>0, \, (J_{12}-J_{23})^2-(1-J_{13})^2>0,\quad \textrm{for $\beta_2=0$}, \\
&(J_{13}+J_{23})^2-(1+J_{12})^2>0, \, (J_{13}-J_{23})^2-(1-J_{12})^2>0,\quad \textrm{for $\beta_3=0$}.
\end{aligned}
\end{equation}
\end{thm}

From Proposition \ref{th:noreturnp} we know that the point of contact which came from an interior of one face of the octahedron and reached one of its edges cannot return to the interior of the face directly from that edge. An important question is how many faces of the octahedron the point of contact can reach in normalized cases?

\begin{thm}\label{prop:threefaces} The point of tangency of the ellipsoid and the octahedron for a lasso solution with normalized data cannot reach the interior of more than two different faces of the octahedron.
\end{thm}
\begin{proof}
Let us suppose the opposite, that it is possible for the point of tangency to reach the interior of three different faces. In order to simplify the exposition and without loss of generality, let us suppose that $\hat{\delta}_i^j=1$, for $j=1,2,3$, and that the point of contact passes over the edges belonging to the planes $\beta_1=0$ and $\beta_2=0$. Let us denote $J_{12}=a, J_{13}=b, J_{23}=c$. From \eqref{eq:usloviprelaz}, the following conditions should be satisfied:
\begin{equation}\label{eq:prelp=3}
-1<a,b,c<1,(a+b+c+1)(a+b-c-1)>0,(a+b-c-1)(a-b-c+1)>0.
\end{equation}
By adding two last relations, one gets $2(a+b-c-1)(a+1)>0$. Since $a+1>0$, one gets
\begin{equation}\label{uslov1}
a+b-c-1>0.
\end{equation}
On the other hand, $s_1$ is the solution of the equation \eqref{eq:solbeta3} with $\beta_1=0$ or
$
-\hat{\beta}_1\Delta=s_1(K_1^{11}+K_1^{12}+K_1^{13}),
$
where $\Delta=\det K_1$ is given by \eqref{delta} and $K_1^{ij}$ is the $(i,j)$-cofactor of the matrix $K_1$. Thus,
the numerator of $s_1$ is $-\hat{\beta}_1\Delta$. From the positive-definiteness of $K_1$, it follows that $\Delta>0$. Since $s_1<0$, and since we supposed that $\hat{\beta_1}>0$, one gets that the denominator of $s_1$ is positive,
$
1-c^2+bc-a+ac-b=(1-c)(-a-b+c+1)>0,
$
or equivalently, $(a+b-c-1)<0$. This leads to the contradiction with \eqref{uslov1}.
\end{proof}

Theorem \ref{thm:sol1} gives the solution of the lasso problem for all $s$ such that $s_I<s<0$. As we mentioned, when $s<s_I$, instead along an interior of a  face, the point of tangency moves along one of the edges $\beta_i=0$, for some $i\in\{1,2,3\}$.
For example, let us suppose that it moves along the edge $\beta_3=0$.
Then, the normal $\bf{n}$ to the ellipsoid is perpendicular to the intersection of the boundary \eqref{eq:region3}  with $\beta_3=0$, that is
\begin{equation}\label{eq:cont2}
|\beta_1|+|\beta_2|=t.
\end{equation}
In other words, the projection of the gradient of the ellipsoid to the plane $\beta_3=0$ is collinear with the vector of the normal to \eqref{eq:cont2} and the coefficient of the proportionality is negative. The projection of the gradient of the ellipsoid for $\beta_3=0$ is
$
(\beta_1+J_{12}\beta_2-b_1,J_{12}\beta_1+\beta_2-b_2).
$

The  condition that the projection is perpendicular  to \eqref{eq:cont2} can be written in the form
\begin{equation}\label{eq:sol3}
\left(\begin{matrix}1&J_{12}\\
        J_{12}&1\end{matrix}\right)\left(\begin{matrix}\beta_1\\ \beta_2\end{matrix}\right)=\left(\begin{matrix}b_1+s\hat{\delta}_i^1\\ b_2+s\hat{\delta}_i^2\end{matrix}\right).
\end{equation}
In the case that $\beta_1=0$, the same procedure gives us
\begin{equation}\label{eq:sol1}
\left(\begin{matrix}1&J_{23}\\
        J_{23}&1\end{matrix}\right)\left(\begin{matrix}\beta_2\\ \beta_3\end{matrix}\right)=\left(\begin{matrix}b_2+s\hat{\delta}_i^2\\ b_3+s\hat{\delta}_i^3\end{matrix}\right).
\end{equation}
If $\beta_2=0$, one gets
\begin{equation}\label{eq:sol2}
\left(\begin{matrix}1&J_{13}\\
        J_{13}&1\end{matrix}\right)\left(\begin{matrix}\beta_1\\ \beta_3\end{matrix}\right)=\left(\begin{matrix}b_1+s\hat{\delta}_i^1\\ b_3+s\hat{\delta}_i^3\end{matrix}\right).
\end{equation}

Summarizing, we get:

\begin{thm}\label{th:edge} For $s<s_I$, when the ellipsoid touches the boundary along the edge $\beta_i=0$, the solution of the lasso reduces to the solution of the two-dimensional lasso problem in the plane $\beta_i=0$ for $i=1$ \eqref{eq:sol1}, for $i=2$ \eqref{eq:sol2}, and for $i=3$ \eqref{eq:sol3}. The corresponding two-dimensional covariance matrix  is obtained from the three-dimensional covariance matrix  \eqref{eq:corrmat3} by omitting  the $i$-th row and the $i$-th column.
\end{thm}

We will show in Example \ref{exam:transfer3d} that in the case $p=3$  in a normalized case, it is possible for the point of tangency of the ellipsoid and the octahedron to pass over one edge, from the interior of one face to the interior of the adjacent face. In $p=2$ case, it was mentioned in \cite{Tib} that it would not be possible. Here we  present an analytic  proof of that fact.
\begin{prop}\label{p2pass} In the case $p=2$,  in the normalized case, the point of tangency of the ellipse and the square cannot belong to the interior of more than one edge of the square.
\end{prop}
\begin{proof} In the case $p=2$,
before the point of tangency reaches a vertex of the square, the solution of the lasso satisfies the linear system of equations \eqref{eq:sol3}. Let us suppose that the point of tangency touch the vertex $\beta_1=0$ for $s=s_1$. In order to pass to the interior of the adjacent edge for the value $s=s_1'$, the following conditions should be satisfied:
$$
s_1=\frac{-\hat{\beta}_{1}(1-J_{12}^2)}{\hat{\delta}_i^1-J_{12}\hat{\delta}_i^2},\quad
s_1'=\frac{-\hat{\beta}_{1}(1-J_{12}^2)}{-\hat{\delta}_i^1-J_{12}\hat{\delta}_i^2},\quad
s_1s_1'>0.
$$
Since the numerators of $s_1$ and $s_1'$ are the same, the product of their denominators should be positive. Hence $J_{12}^2-1>0$. This leads to the contradiction with the fact $|J_{12}|<1$.
\end{proof}

Proposition \ref{p2pass} uses the fact that the semiaxes of the ellipses are parallel to the edges of the square. This is a consequence of the assumption that the matrix $K_1$ is normalized: $J_{11}=J_{22}=1$.  The proposition is not valid without the normalization assumption, i.e. in the case, when $J_{11}$ and $J_{22}$ are not both equal to $1$, as it follows from \cite{MY}.

Let us assume now that $K_1$ is not normalized, but it is still positive-definite:
\begin{equation}\label{k1nonnorm}
K_1=\left(\begin{matrix}J_{11}&J_{12}\\
                        J_{12}&J_{22} \end{matrix}\right).
\end{equation}
The positive-definitness gives
\begin{equation}\label{eq:delta2}
J_{11}>0,\quad \Delta_2=J_{11}J_{22}-J_{12}^2>0.
\end{equation}
Let $J_{13}, J_{23}, J_{33}$ be given such that $3\times3$ matrix $K=(J_{ij})$ is positive definite.
Let us also suppose that there are ellipses from the family, which are tangent to two edges of the squares that belong to two different quadrants.
According to previous considerations, there are two rays, both initiated at the center of the ellipses, such that for $s_I<s<0$ the solutions of the lasso belong to the first ray, while for $s<s_I'<s_I<0$ the solutions belong to the second ray, for some values $s_I$, $s_I'$. Let us denote by $(\xi, 0)$ and $(0,\eta)$ the points where these rays intersect the corresponding coordinate lines.
Without loss of generality, let us suppose
\begin{equation}\label{uslovip2prelaz}
\hat{\beta}_1>0,\quad\hat{\beta}_2>0,\quad\xi>0\quad\eta<0\quad |\xi|>|\eta|,
\end{equation}
see Figure  \ref{pic:2dnonnormalized}.
Using $K_1(\beta-\hat{\beta})=s\delta$, and $K_1^{-1}(J_{13}\ J_{23})^T=\hat{\beta}$ one calculates
$$
\hat{\beta}_1=\frac{J_{22}J_{13}-J_{12}J_{23}}{\Delta_2},\quad \hat{\beta}_2=\frac{-J_{12}J_{13}+J_{22}J_{23}}{\Delta_2}, \quad \xi=\frac{J_{13}-J_{23}}{J_{11}-J_{12}},\quad \eta=\frac{J_{13}+J_{23}}{J_{12}+J_{22}},
$$
where $\Delta_2$ is the determinant of $K_1$, with an explicit expression given in \eqref{eq:delta2}.
Then, the conditions \eqref{uslovip2prelaz} can be rewritten as:
\begin{equation}\label{uslp2novi}
\begin{aligned}
J_{22}J_{13}-J_{12}J_{23}&>0,\quad -J_{12}J_{13}+J_{22}J_{23}>0,\quad\frac{J_{13}-J_{23}}{J_{11}-J_{12}}>0,\\
\frac{J_{13}+J_{23}}{J_{12}+J_{22}}&<0,\quad\Big|\frac{J_{13}-J_{23}}{J_{11}-J_{12}}\Big|>\Big|\frac{J_{13}+J_{23}}{J_{12}+J_{22}}\Big|.
\end{aligned}
\end{equation}

Now, we can formulate the following

\begin{prop}\label{prop:2nonnorm}  In the case $p=2$, with the matrix $K_1$ given by \eqref{k1nonnorm}, which is not normalized (i.e. assuming  $J_{11}\ne1$ or $J_{22}\ne1$), suppose that conditions \eqref{uslp2novi} are satisfied. If the following conditions
\begin{equation}\label{uslprl2}
\frac{J_{22}J_{13}-J_{12}J_{23}}{\Delta_2}>\frac{J_{13}-J_{23}}{J_{11}-J_{12}},
\frac{J_{22}J_{13}-J_{12}J_{23}}{\Delta_2}>\frac{-J_{12}J_{13}+J_{22}J_{23}}{\Delta_2}+\frac{J_{13}-J_{23}}{J_{11}-J_{12}},
\end{equation}
are satisfied, then a passing over a vertex happens, otherwise, a passing does not occur.
\end{prop}

\begin{proof}
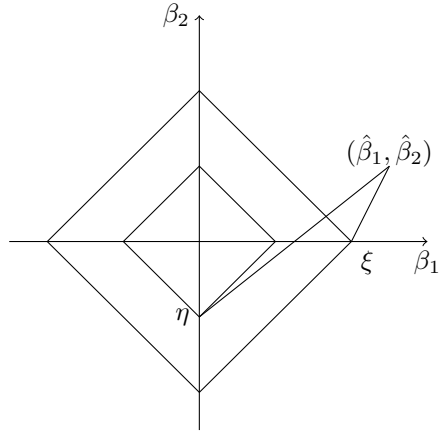
\begin{figure}[h]\centering
\begin{tikzpicture}
\draw[->] (-2.5,0) -- (3,0) node[anchor= north ] {$\beta_1$ };
\draw[->] (0,-2.5) -- (0,3) node[anchor= east] {$\beta_2$};
\draw (-1,0)--(0,-1)--(1,0)--(0,1)--(-1,0);
\draw (-2,0)--(0,-2)--(2,0)--(0,2)--(-2,0);
\draw (2.5,1)--(2,0) node[anchor= north west] {$\xi$};
\draw (2.5,1)--(0,-1) node[anchor= east] {$\eta$};
\path (2.5,1.2) node {$(\hat{\beta}_1,\hat{\beta}_2)$};
\end{tikzpicture}
\caption{ The case $p=2$} \label{pic:2dnonnormalized}
\end{figure}

Let us suppose that $\xi\ne\hat{\beta}_1, \eta\ne\hat{\beta}_2$. The parameter $s$ that corresponds to the point of intersection of the $O\beta_1$ axis and the ray that connects the center with the point $(0,\eta)$, is smaller than the parameter $s$ that corresponds to the point $(\xi,0)$. Thus, in the case $\hat{\beta}_1>\xi$, the slope of the ray that connects the center with $(\xi,0)$ is greater  than the slope of the line that connects points $(\xi,0)$ and $(0,\eta)$.

One has $\frac{\hat{\beta}_2}{\hat{\beta}_1-\xi}>\frac{-\eta}{\xi}$. Summarizing, one gets that
\begin{equation}\label{eq:p=2cond}
(\hat\beta_1-\xi)\Big(\frac{\eta}{\xi}+\frac{\hat\beta_2}{\hat\beta_1-\xi}\Big)>0.
\end{equation}
From the conditions of the tangency of the ellipses with the edges of the corresponding squares at the points $(\xi,0)$ and $(0,\eta)$, one gets the system of equations
that can be written in the form
$$
\frac{J_{22}}{J_{11}}\hat\beta_2-\frac{J_{12}}{J_{11}}(\xi-\hat\beta_1+\hat\beta_2)=\hat\beta_1-\xi,
\frac{J_{22}}{J_{11}}(\eta-\hat\beta_2)+\frac{J_{12}}{J_{11}}(\eta-\hat\beta_1-\hat\beta_2)=\hat\beta_1.
$$
The solution of this system is
$$
\frac{J_{22}}{J_{11}}=\frac{(\hat\beta_1-\xi)(\eta-2\hat\beta_1)+\xi\hat\beta_2}{(\eta-\hat\beta_2)(2\hat\beta_2+\xi)-\hat\beta_1\eta},\,
\frac{J_{12}}{J_{11}}=\frac{\hat\beta_2\hat\beta_1-(\hat\beta_1-\xi)(\eta-\hat\beta_2)}{(\eta-\hat\beta_2)(2\hat\beta_2+\xi)-\hat\beta_1\eta}.
$$
Since the matrix $K_1$ is positive-definite, the condition $\frac{J_{22}}{J_{11}}-\Big(\frac{J_{12}}{J_{11}}\Big)^2>0$ should be satisfied.
One calculates that
$$
\frac{J_{22}}{J_{11}}-\Big(\frac{J_{12}}{J_{11}}\Big)^2=\frac{2}{D^2}(\eta-\hat\beta_2-\hat\beta_1)(\hat\beta_2-\hat\beta_1+
\xi)(\eta\hat\beta_1-\xi\eta+\xi\hat\beta_2),
$$
where $D=(\eta-\hat\beta_2)(2\hat\beta_2+\xi)-\hat\beta_1\eta$. When  $\hat{\beta}_1>\xi$,
 the expression given in the first bracket is negative.
One has  $\eta\hat\beta_1-\xi\eta+\xi\hat\beta_2=\xi(\hat\beta_1-\xi)\big(\frac{\eta}{\xi}+\frac{\hat\beta_2}{\hat\beta_1-\xi}\Big)$. Since $\xi>0$ using \eqref{eq:p=2cond}, we have that the third bracket is positive. Since,
 $K_1$ is positive-definite, the passing may happen if the second bracket $\hat\beta_2-\hat\beta_1+\xi$ is negative. These conditions are equivalent to \eqref{uslprl2}. Otherwise, one gets the contradiction. Let us mention that the last expression is negative when the slope of the line that connects the center with the point $(\xi,0)$ is less than $1$.

When $\hat{\beta}_1<\xi$, using that $\eta<0$ and the obvious relation $\eta\hat\beta_1-\xi\eta+\xi\hat\beta_2=\eta(\hat\beta_1-\xi)+\xi\hat\beta_2$ one concludes that the first bracket is negative, while the second and the third brackets are positive. One gets the contradiction.

 When $\hat{\beta}_1=\xi$, one gets the solutions
$$
J_{22}=J_{12},\quad J_{22}=\frac{J_{11}\hat{\beta}_1}{2\eta-2\hat{\beta}_2-\hat{\beta}_1}.
$$
From the positive-definiteness, one gets $\hat{\beta}_1+\hat{\beta}_2<\eta$, that is not the case. The case $\hat{\beta}_2=\eta$ can be treated similarly.

\end{proof}

\begin{exm}\label{exam:transfer3d} Now we will present an example  of a normalized data for $p=3$ in which the point of tangency of the ellipsoid and the octahedron can pass from  the interior of one face of the octahedron to  the interior of an adjacent face.  This will exemplify
a resurgence of the coefficient $\beta_1$  in an explicit form for the first time  for normalized data (to the best of our knowledge), and thus, provides a proof that such a phenomenon is indeed possible.

 Let the hypothetical data  are given:
$$
X^T=\left(\begin{matrix}  -1.62553679& -1.45872292& -1.29498175\\
 0.81876630& 1.12398462& -1.34295625\\
 0.43204938 &-0.32403703& 0.10801234\\
 0.02419136 &-0.06153051& -0.08845401\\
 1.62553679 &1.45872292 &1.29498175\\
 -0.81876630 &-1.12398462 &1.34295625\\
 -0.43204938 &0.32403703 &-0.10801234\\
 -0.02419136 &0.06153051 &0.08845401
 \end{matrix}\right),
$$
with the vector of the responses
\begin{multline*}
Y=(-1.71048007, -0.71992247, -0.21602469, 0.09645148,\\
 1.71048007, 0.71992247,  0.21602469, -0.09645148)^T.
\end{multline*}
Then, the covariance matrix  of the predictors and its inverse are
$$
K_1=\left(\begin{matrix}
1&0.9&0.3\\
0.9&1&0.1\\
0.3&0.1&1
\end{matrix}
\right),\quad
K_1^{-1}=\left(\begin{matrix}
6.875&-6.041667&-1.45833\\
-6.041667&6.31944&1.180556\\
-1.45833&1.180556&1.31944
\end{matrix}
\right),
$$
and $J_{14}=0.6$, $J_{24}=0.5$, $J_{34}=0.9$, $J_{44}=1$.

The center of the ellipsoids \eqref{coefeq} is at the point $\hat{\beta}=(-0.20833, 0.59722, 0.90278)$. Hence $\delta_i^1=-1$ while $\delta_i^2=\delta_i^3=1$. The values $s_1$, $s_2$, $s_3$ are calculated from \eqref{eq:solbeta3} as the conditions that $\beta_1=0$, $\beta_2=0$ and $\beta_3=0$ respectively. One gets $s_1= -0.014493$, $s_2=-0.044102$, $s_3=-0.228070175$. Thus, for $-0.014493<s<0$, the lasso solution is
$
\beta_1=-0.20833-14.375s,\, \beta_2=0.59722+13.541667s,\,  \beta_3=0.902778+3.95833s.
$
For $s=-0.014493$, the coefficient $\beta_1=0$. For $s$ less and close enough to $-0.014493$, the contact point of the ellipsoid and the octahedron is on the edge $\beta_1=0$. The question is which of the two scenarios will happen: will the point of contact reach a vertex or will it pass to the adjacent face? Each of three criteria \eqref{eq:cond1}, \eqref{eq:cond2}, \eqref{eq:cond3} is satisfied here. From \eqref{s1prim}, one gets that the value of $s$, for which the  passing of the contact point from the edge to the interior of the adjacent face happens for $s_1'=-0.3333$. On the other hand, for $s<-0.014493$ the solution of the lasso is obtained from an induced lasso for the reduced dimension $p=2$.
The covariance matrix $
K_1^{(2)}$ of the induced $p=2$ lasso is $
{K_1^{(2)}}_{11}=
{K_1^{(2)}}_{22}=1$ and $
{K_1^{(2)}}_{12}=
{K_1^{(2)}}_{21}=0.1$

Using \eqref{eq:sol1} one gets
\begin{equation}\label{hodpoivici}
\beta_1=0,\quad\beta_2=0.41414+0.90909s, \quad \beta_3=0.85859+0.90909s.
\end{equation}
The values for $s$ for which $\beta_2=0$ and $\beta_3=0$ are $s^{(2)}_2=-0.454546$ and $s^{(2)}_3=-0.9444$ respectively. Since  $s_1'>s^{(2)}_2>s^{(2)}_3$, the transition of the contact point from the edge to the adjacent face will indeed happen here. One gets that for all $s$ such that  $-0.3333<s<-0.014493$, the solution to the lasso is \eqref{hodpoivici}, and for $s=-0.3333$, the ellipsoid will be tangent to the adjacent face of the octahedron that shares the edge $\beta_1=0$ with the initial face. Thus, for $s$ less and close enough to $-0.3333$,
the point of contact will belong  to the interior the face of the octahedron which belongs to the first octant. There $\beta_1$, which was equal to zero along the edge, is again different from zero, but now it is positive. Since the point of contact is now in the first octant, all delta's are now equal to +1. The values $s^{(3)}_2$, $s^{(3)}_3$ are calculated from \eqref{eq:solbeta3} and the conditions that $\beta_2=0$ and $\beta_3=0$ respectively.
One gets $s^{(3)}_2=-0.4095$, $s^{(3)}_3=-0.8667$. Thus, for $-0.4095<s<-0.3333$, the lasso solution is
$
\beta_1=-0.20833-0.625s,\, \beta_2=0.59722+1.45833s,\,  \beta_3=0.902778+1.041667s.
$
When $s=-0.4095$, the coefficient $\beta_2=0$.  Now, the condition \eqref{eq:cond2} for the transition of the point of contact from the edge $\beta_2=0$ to the interior of the adjacent face is not satisfied. For $s<-0.4095$, the point of contact that represents the solution will
be on the edge $\beta_2=0$, until it hits a vertex of the octahedron. The solution of the lasso reduces again to an induced lasso problem of the reduced  dimension $p=2$. Using \eqref{eq:sol2}, for $-0.47143<s<-0.4095$, the lasso solution is
$
\beta_1=0.362637+0.76923s,\, \beta_2=0,\, \beta_3=0.791209+0.76923s.
$
Finally, for all $s$ such that $s<-0.47143$, the solutions are $\beta_1=0$, $\beta_2=0$.

We have also calculated this example using the R. The comparison of the results obtained above with the results calculated using the R are given in the following table. In order to compare the values of $s$ in our formulas with the values of $\lambda$ in the results obtained using the R software, we first fix $\beta_1$ to be the same in both cases. Then we calculate the remaining coefficients $\beta$ in two ways. We see an excellent agreement. The ratio $-s/\lambda$ is approximatively equal to $8$, which is equal to $N$ in this case.
\begin{table}[h]
\caption{Comparison between the solutions obtained using our formulas and solutions obtained using R-software}
\begin{tabular}{|c|c|c|c|c|}
\hline
\multicolumn{4}{|c|}{Formulas}&\\
\hline
$s$&$\beta_1$&$\beta_2$&$\beta_3$&\\
\hline
 -0.002907005&-0.166541804&0.557854307&0.891271115&\\
\hline
-0.007237362&-0.10429292&0.499214053&0.874130133&\\
\hline
-0.011456138&-0.043648014&0.442084792&0.857430825&\\
\hline
-0.028792155&0&0.38796534&0.83241534&\\
\hline
-0.087929435&0&0.33420423&0.77865423&\\
\hline
-0.3920808&0.0367205&0.025436807&0.494360369&\\
\hline\hline
\multicolumn{4}{|c|}{R-software}&\\
\hline
$\lambda$&$\beta_1$&$\beta_2$&$\beta_3$&$s/\lambda$\\
\hline
0.0003531643&-0.166541804&0.557949387&0.891354458 &-8.23131\\
\hline
0.000891543&-0.10429292&0.49932052&0.87422348&-8.11779\\
\hline
0.001419589&-0.043648014&0.442185553&0.857519137&-8.07004\\
\hline
0.04869856&0&0.38796534&0.83241006&-8.05118\\
\hline
0.003599176&0&0.33420423&0.77864868&-7.99965\\
\hline
0.01099136&0.03672050&0.02771257&0.496624408&-7.99987\\
\hline
\end{tabular}
\end{table}
One can see that by decreasing $s$ (or increasing $\lambda$) the coefficient $\beta_1$  goes from a negative value to zero, and after a while being zero, again starts growing, and gets positive values. On the  Figure \ref{sl:4.5R} the coefficients are presented as functions of $\log \lambda$. The phenomenon of passing from one octant to another one can be clearly observed.
\begin{figure}[h]
\includegraphics[width=8cm,height=5.5cm]{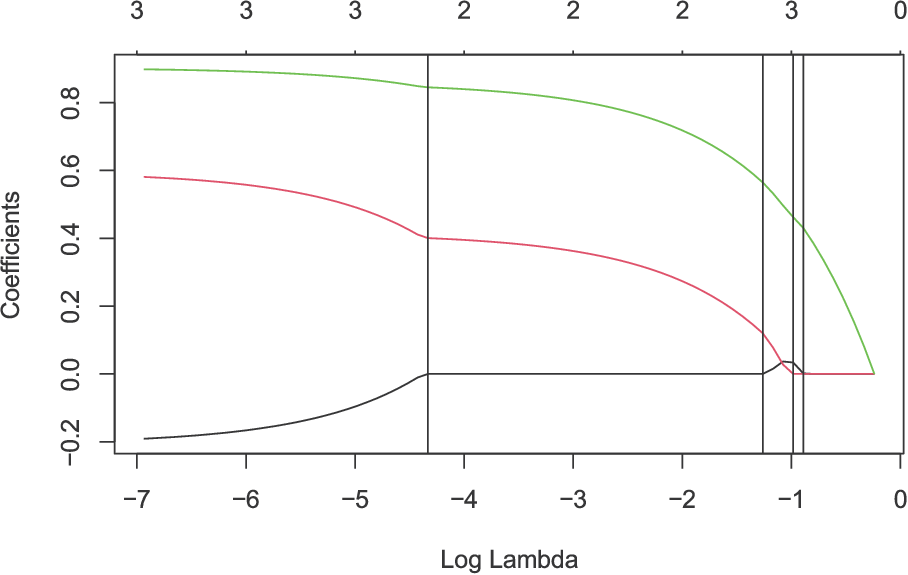}
\caption{Coefficients $\beta$ as functions of $\log \lambda$. The resurgence of $\beta_1$ is manifested through a ``bump" around $-1$.}
\label{sl:4.5R}

\end{figure}
\end{exm}

We illustrate Theorem \ref{thm:raysfinal} using  \textsc{Example} \ref{exam:transfer3d}.
In Figure \ref{picexm4.4a} we see two rays form the center $\hat{\beta}$, since in this example, there is a passing over one of the edges (the edge $\beta_1=0$).

\begin{figure}[h] \centering
\usetikzlibrary{3d}
\tdplotsetmaincoords{70}{130}
\begin{tikzpicture}[scale=4.5, tdplot_main_coords]\label{picexm4.4a}
\path (-0.2,0.6,0.96) node[]{\footnotesize $(\hat\beta_1,\hat\beta_2,\hat\beta_3)$};
\draw[very thick,->](-0.20833,0.59722,0.902778)--(0, 0.4009615, 0.8454099);
\draw[very thick,->](0, 0.4009615, 0.8454099)--(0,0.111413,0.555863);
\draw[very thick,->](0,0.111413,0.555863)--(0.0476075,0,0.47621536);
\draw[very thick,->](0.0476075,0,0.47621536)--(0,0,0.428902);
\draw[very thick,->](0,0,0.428902)--(0,0,0);
\draw[->] (-0.5,0,0)->(0.5,0,0) node[anchor=north east]{$\beta_1$};
\draw[->] (0,0,0)->(0,0.7,0) node[anchor=north west] {$\beta_2$};
\draw[->] (0,0,0)->(0,0,0.7) node[anchor=south]{$\beta_3$};
\draw[dashed] (0,0.111413,0.555863) -- (-0.20833,0.59722,0.902778);
\filldraw (-0.20833,0.59722,0.902778) circle (0.3pt);
\filldraw (0, 0.4009615, 0.8454099) circle (0.3pt) ;
\filldraw (0,0.111413,0.555863) circle (0.3pt) ;
\filldraw (0.0476075,0,0.47621536) circle (0.3pt) ;
\filldraw (0,0,0.428902) circle (0.3pt) ;
\filldraw (0,0,0) circle (0.3pt) ;

\end{tikzpicture}

\caption{\textsc{The solution of lasso from Example} \ref{exam:transfer3d}}

\end{figure}

\begin{remark}

From the point of view of a lasso procedure, one can start with  a covariance or correlation matrix. However, there are infinitely many data sets with the same statistic properties (see for example \cite{MF}),  i.e. those which produce the same covariance matrix. Since the lasso implementation in the R packages assumes the data as an input, rather than a covariance matrix, let us explain briefly one of the many possible ways to find a centered data set that has the given matrix of covariances $K_1$.  Given a normalized positive-definite matrix $K_1$,  one needs to find a centered matrix $X$, such that $XX^T=K_1$. Since $K_1$ is positive definite matrix, there is an orthogonal matrix $S$, such that $K_1=SDS^T$, where $D$ is diagonal matrix with the eigenvalues of $K_1$ on diagonal. Using positive definiteness of $K_1$, one gets $K_1=S\sqrt{D}S^TS\sqrt{D}S^T$. Thus, a candidate for $X$ is the matrix $\tilde{X}=S\sqrt{D}S^T$. Since $K_1$ is normalized, its diagonal elements are equal to $1$, and the sum of the squares of the elements of each row of $\tilde{X}$ is equal to one. But, is may happen that the mean of each row is not equal to zero. Then, a centered $X$ can be obtain as a $p\times2p$ block matrix $X=[\frac{1}{\sqrt{2}}\tilde{X}|-\frac{1}{\sqrt{2}}\tilde{X}]$. It is easy to check that the rows of $X$ have the first two moments equal to $0$ and $1$ respectively.

Instead of the spectral decomposition used above, for example one may use the Cholesky decomposition  (see e.g. \cite{HJ, Ge}) in order to get another matrix $X$, such that $XX^T=K_1$ and its rows have the first two moments equal to $0$ and $1$ respectively.
\end{remark}

We conclude our study of the case $p=3$ by giving necessary and sufficient conditions for the point of tangency of the ellipsoids and the octahedrons to pass from the interior of one face to the interior of an adjacent face. We will assume that the common edge of the two faces belongs to the intersection of the plane $\beta_1=0$ with the octahedron. One can easily rewrite the conditions of the following Proposition for the cases of passings over the edges which belong to the plane $\beta_2=0$ or the plane $\beta_3=0$.
\begin{prop}\label{prop:transfer3d} Let $s_1$, $s_2$ and $s_3$ be given by formulas \eqref{eq:sfor3}. The point of tangency of the ellipsoid and the octahedron will pass over an edge belonging to the plane $\beta_1=0$  in a normalized case if and only if the following conditions are satisfied:
\begin{enumerate}
\item\label{pid1} $s_1<0$ and one of the following conditions is satisfied:
\begin{enumerate}
\item $s_2, s_3<0$, and $s_2<s_1$, $s_3<s_1$;
\item $s_3>0$, and $s_2<s_1<0$;
\item $s_2>0$, and $s_3<s_1<0$;
\item $s_2,s_3$ are both positive;
\end{enumerate}
\item\label{pid2} $s_1'<0$, where $s_1'$ is given by \eqref{s1prim};
\item\label{pid3} $s^{(2)}<s_1'$, where $s^{(2)}$ is the one of $s_2^{(2)}$, $s_3^{(2)}$ that is maximal negative, where
$$
s_2^{(2)}=\frac{-J_{24}+J_{23}J_{34}}{\hat{\delta}_i^2-\hat{\delta}_i^3J_{23}},\quad  s_3^{(2)}=\frac{J_{34}+J_{23}J_{24}}{\hat{\delta}_i^3-\hat{\delta}_i^2J_{23}}.
$$
\end{enumerate}
\end{prop}
\begin{proof} We already proved that \eqref{pid1} and \eqref{pid2} are necessary conditions for passing across an edge. The point of tangency  passes to the interior of an adjacent face if and only if the value of $s_1'$ is greater than the values of $s$ that correspond to the vertices $\beta_1=\beta_2=0$ and $\beta_1=\beta_3=0$. Using Theorem \ref{th:edge} and \eqref{eq:sol1}, one gets the values $s_2^{(2)}$ and $s_3^{(2)}$. Let the conditions \eqref{pid1}, \eqref{pid2} and \eqref{pid3} be satisfied. The condition \eqref{pid1} gives that the point of tangency  will reach the edge $\beta_1=0$ first, among all the edges of the face. The conditions \eqref{pid2} mean that a passing to the interior of the adjacent face is possible, and it will occur if and only if the condition \eqref{pid3} is satisfied.
\end{proof}

If a data set is nonnormalized, in the case $p=2$, we provided an explicit example that admits passing over a vertex. In \cite{MY}, an example was provided that for $p=3$ had three passings across the edges.

\section{ Complexity in the maximal dimension of lasso for general $p$}\label{sec:compmax}

Having in mind the Proposition \ref{th:noreturnp}, a natural question is: What is the maximal number  of segments that a simple polynomial chain that represents the lasso solution may have? A closely related question is: how many interiors of orthants of the maximal dimension $p$ can a lasso solution intersect?

Let us recall that the polytope defined by the $\ell_1$ norm $B(t)$ is a cross-polytope \cite{Co, Zi}. We code the path of a lasso solution with a sequence of $p$-dimensional vectors $\delta_i$ with $1,-1,0$ as entries. Here $\delta_i$ is, as before, the vector of the signs of a segment of the solution.

The entry $\delta_i^j=0$ means that $\beta_j=0$ for the $i$-th segment. The exception is $\delta_i=0$, when the segment degenerates to the origin. Since there are three possible choices for $\delta_i^j$, it follows that the maximal number of $\delta_i$'s is not greater than $3^p$. However, in \cite{MY}, it was proved analytically that, if the vector $\delta_i$ appears in the solution, the vector $-\delta_i$ does not appear, hence the maximal number of the segments in the lasso solution cannot exceed $(3^p+1)/2$. It can be easily seen from the geometric point of view that appearing of both vectors  $\delta_i$ and $-\delta_i$ would mean that there are two homothetic ellipsoids from the family that are externally tangent to two parallel faces of two cross polytopes, and these faces belong to two different orthants,  leading to the contradiction. As it was shown in \cite{MY}, the estimate from above of $(3^p+1)/2$ segments of the lasso solution is sharp: a remarkable construction of a lasso with the solution that has exactly $(3^p+1)/2$ segments in an arbitrary dimension $p$ was given in \cite{MY}. This also gives a sharp estimate on the number of segments that belong (except the end point) to the interior of orthants of the maximal dimension $p$ is $2^{p-1}$. Every such segment is coded with a vector $\delta$ that does not contain any entry equal to $0$.

A direct consequence of the fact that two opposite $\delta$'s cannot both appear in the solution is the following statement.

\begin{prop}{\label{prop:commonedge}} If two segments of a solution belong to two orthants of the maximal dimension $p$, then the dimension of the intersection of the closures of these orthants is at least one.
\end{prop}

\begin{cor}\label{cor:5.1} For arbitrary pair of vectors $\delta_i$, $\delta_j$ without zero entries in the sequence representing a lasso solution, there is a spot, say $k$, such that $\delta_i^k=\delta_j^k$.
\end{cor}

 In what follows, we will present an upper bound for the number of the segments intersecting the interior of orthants of the maximal dimension $p$ of the lasso solutions in the normalized case for any $p$.

\begin{prop}\label{prop:edgenorm} In the normalized case in an arbitrary dimension $p$, if a segment of a lasso solution belongs to one of the coordinate axis, then it contains the origin. Thus, such a segment is unique in a solution, if it exists.
\end{prop}

For the existence question see Remark \ref{rem:degen} below.

\begin{proof} Let us suppose opposite, that a segment with the endpoints $A$ and $B$ different from the origin appears in the lasso solution, and that $AB$ belongs to a coordinate axis, say $O\beta_i$ axis. Then there are two segments of the solution $AC$ and $BD$, such that $AC$ belongs to a two-dimensional plane $\beta_iO\beta_j$ and $BD$ belongs to a plane $\beta_iO\beta_k$ for some $j$ and $k$. This means that a passing over an edge of the cross-polytope happens. This gives the contradiction, since in normalized cases, such a crossing  never happens.
\end{proof}

\begin{remark}\label{rem:degen} It can happen that two (or more) coefficients of $\beta_i$ become zero for the same value of the parameter $s$. In such a case, it is possible  that a segment that belongs to a two-dimensional coordinate plane has an endpoint at the origin, and the segment of the lasso solution that belongs to a coordinate axis then does not appear in the solution of a lasso. However, in a situation of a general position, this does not occur.
\end{remark}

Proposition \ref{prop:edgenorm} reflects a key difference between lasso solutions for normalized vs. non-normalized data: in the normalized case there is no passing over an edge of the cross-polytope $B(t)$.

\begin{prop}\label{prop:intersection2} Let us suppose that two segments of a solution of a lasso, with a general data, belong to two orthants of the maximal dimension $p$, such that their closures intersect over a ray $O\beta_i$ for some $1\leq i\leq p$. Then, there exists a point $P$ that belongs to the intersection of the lasso path and the coordinate line $O\beta_i$.
 \end{prop}
\begin{proof} Without loss of generality, we consider the case of two segments of a solution with their assigned vectors $\delta_1=(1,1,...,1)^T$ and $\delta_2=(1,-1,-1,...,-1)^T$. There are parameters $s_1$ and $s_2$ and the corresponding solutions vectors $\beta$ and $\eta$  such that
$
K_1(\beta-\hat{\beta})=s_1\delta_1,\quad K_1(\eta-\hat{\beta})=s_2\delta_2.
$
By multiplying the first equation by $s_2/(s_1+s_2)$ and the second one by $s_1/(s_1+s_2)$, after adding them, one gets
$
K_1(\xi-\hat{\beta})=s_3\delta,
$
 where $s_3=\frac{2s_1s_2}{s_1+s_2}$, $\delta=(1,0,0,...,0)^T$ and $\xi=\frac{s_2\beta+s_1\eta}{s_1+s_2}$. Since $s_3$ is between $s_1$ and $s_2$, it follows that there is a segment of the solution of the lasso that belongs to the coordinate line $O\beta_1$. For a one-dimensional lasso the ellipsoid reduces to a segment with the center $b_1$.  We define $\hat{\xi}$ as
  $$
  1\cdot(\xi_1-\beta_1)+\dots+(K_1)_{1p}(\xi_p-\beta_p)=\hat{\xi}-b_1=s_3.
  $$
It gives the solution of the one-dimensional lasso. Since $s_3$ is between $s_1$ and $s_2$, the endpoints of this segment are different from the origin, giving  the contradiction with  Proposition \ref{prop:edgenorm}.
\end{proof}

\begin{thm}\label{thm:2orthants} Let us consider a lasso with a normalized data: $K_{1,ii}=1$. If two segments of the solution belong to two orthants of the maximal dimension $p$, then the dimension of the intersection of the closures of these orthants is at least two.
\end{thm}
\begin{proof} Proposition  \ref{prop:commonedge} gives that the dimension of the intersection of the closures of these orthants is at least one. If it is one, then from Proposition \ref{prop:intersection2}, it follows that the path of the solution contains a segment that belongs to a coordinate axis. In other words, a passing over an edge of the cross-polytope $B(t)$ happens. This is not possible in a normalized case. Thus, the dimension of the intersection of two $p$-dimensional orthants that contain segments of a solution is at least two.
\end{proof}
In terms of vectors $\delta$, Theorem \ref{thm:2orthants} has the following interpretation

\begin{cor}\label{cor:5.2} In a normalized case, if vectors $\delta_i$ and $\delta_j$ with all coordinates different from zero appear in a solution, then at least two of their coordinates coincide.
 \end{cor}

\section{Combinatorics of  binary codes and graphs}\label{sec:comb}
 {

Denote by $h_{p,r}$ the maximal number of binary words of length $p$ such that every two words match at least at $r$ spots, $r\le p$.

 In terms of binary codes, the problem of finding an upper bound for the number of segments that intersect the interiors of orthants of the maximal dimension $p$ of a lasso  solution can be reformulated as follows.

 \begin{prop}\label{prop:upperbp}
 \begin{enumerate}
 \item The number $h_{p,1}$ is a sharp upper bound  for the number of segments that intersect the interiors of orthants of the maximal dimension $p$ of a lasso  solution in nonnormalized cases.
  \item The number $h_{p,2}$ is an upper bound  for the number of segments that intersect the interiors of orthants of the maximal dimension $p$ of a lasso  solution in normalized cases.
 \end{enumerate}
\end{prop}

\begin{proof} Part (1) follows for Corollary \ref{cor:5.1}, while part (2) follows from Corollary \ref{cor:5.2}.
\end{proof}

\begin{remark} For all $p$, $h_{p,1}=2^{p-1}$.
If a word belongs to a maximal set of binary words of length $p$, such that every two words coincide at least at one spot, then its complementary word does not belong to that set. Thus $h_{p,1}\le 2^{p-1}$.
Consider the set of all $p$ words with $1$ at the first spot. Every two words from that set match at least at one spot. There are $2^{p-1}$ words in this set. Thus, $h_{p,1}\ge 2^{p-1}$. This proves the statement.
\end{remark}

\begin{thm}\label{lbound}
\begin{enumerate}
\item
For $p=2\ell$:  $h_{2\ell, 2}\ge2^{2\ell-1}-\frac{1}{2}\binom{2\ell}{\ell}$.
\item
For  $p=2\ell+1$: $h_{2\ell+1, 2}\ge2^{2\ell}-\binom{2\ell}{\ell}$.
\end{enumerate}
\end{thm}
\begin{proof} Assume first $p=2\ell$. We choose the set of $p$-tuples, each of them has at least $\ell+1$ entries equal to $1$. Any two of such words match at least at two spots with $1$ assigned to them. The number of these words is $2^{2\ell-1}-\frac{1}{2}\binom{2\ell}{\ell}$. This proves part (1). When $p=2\ell+1$ is odd, we choose the set $M_{2a}$ of all $p$-tuples with at least $\ell+2$ coordinates equal to $1$. Their number is $2^{2\ell}-\binom{2\ell+1}{\ell+1}$. Any two of these words match on at least two spots. Now, consider the set $M_{2b}$ of all $p$-words with exactly $\ell+1$ coordinates equal to $1$ and with $-1$ at the first spot. Thus, $\ell+1$ of the remaining  $2\ell$ coordinates are equal to $1$ for elements of $M_{2b}$. Any two words from $M_{2b}$ match at least at three spots. It is easy to conclude that any two words, one from $M_{2a}$ and another from $M_{2b}$, match at least at two spots. The number of elements of $M_{2b}$ is $\binom{2\ell}{\ell+1}$. Since $M_{2a}$ and $M_{2b}$ are disjoint, one gets for their union $|M_{2a}|+|M_{2b}|=2^{2\ell}-\binom{2\ell+1}{\ell+1}+\binom{2\ell}{\ell+1}$. Thus, $h_{2\ell+1, 2}\geq 2^{2\ell}-\binom{2\ell}{\ell}$.
\end{proof}

Denote by $e_1, \dots, e_p$ the standard basis of $\mathbb R^p$. We consider the Cayley graphs $G_{2\ell}=\Cay(\mathbb Z_2^{2\ell}, S_{2\ell})$, where $S_{2\ell}=\{e_1,\dots, e_{2\ell}, e_1+\dots + e_{2\ell}\}$ and the Cayley graphs $G_{2\ell+1}=\Cay(\mathbb Z_2^{2\ell+1}, S_{2\ell+1})$, where $S_{2\ell+1}=\{e_1,\dots, e_{2\ell+1}, e_1+\dots + e_{2\ell}\}$. Let us recall that the independence number $\alpha(G)$ of a graph $G$ is the maximal number of vertices such that no pair of those vectors is connected with an edge of the graph $G$. The independence number of a graph $G$ is equal to the clique number of the complement graph $\hat G$, where the clique number is the maximal number of vertices of a complete subgraph of $\hat G$, see \cite{GR}.

\begin{thm}
\begin{enumerate}
\item
For $p=2\ell$:  $h_{2\ell, 2}=\alpha(G_{2\ell})$, where $\alpha(G_{2\ell})$ is the independence number of  the Cayley graphs $G_{2\ell}=\Cay(\mathbb Z_2^{2\ell}, S_{2\ell})$.
\item
For  $p=2\ell+1$: $h_{2\ell+1, 2}=\alpha(G_{2\ell+1})$, where $\alpha(G_{2\ell+1})$ is the independence number of  the Cayley graphs $G_{2\ell+1}=\Cay(\mathbb Z_2^{2\ell+1}, S_{2\ell+1})$.
\end{enumerate}
\end{thm}
\begin{proof}
Let us consider the graph $\hat G_p$ of all the words of length $p$, with two words being connected if and only if they match at least at two spots.
\begin{lem}
\begin{enumerate}
\item
For $p=2\ell$: The graph $\hat G_{2\ell}$ is the complement of the Cayley graphs $G_{2\ell}=\Cay(\mathbb Z_2^{2\ell}, S_{2\ell})$.
\item
For  $p=2\ell+1$: The graph $\hat G_{2\ell+1}$ is the complement of  the Cayley graphs $G_{2\ell+1}=\Cay(\mathbb Z_2^{2\ell+1}, S_{2\ell+1})$.
\end{enumerate}
\end{lem}
\begin{proof} In the Cayley graph $G(\hat{S})=\Cay(\mathbb Z_2^{p}, \hat{S})$ where $\hat{S}=\{e_2+e_3+...+e_p, e_1+e_3+...+e_p,...,e_1+...+e_{p-1},e_1+e_2+...+e_p\}$, two vertices are adjacent if and only if the corresponding words of length $p$ match at exactly one spot, or if they do not match at all. Hence in its complement, two vertices are connected  if and only if they match at least at two spots. Thus, the complement of the graph  $G(\hat{S})$ is $\hat G_{p}$. By changing the basis, we get that for all $p$, the graphs $G_{2p}$ and $G(\hat{S})$ are isomorphic.  This proves the Lemma.

\end{proof}
It is obvious that $h_{p,2}$ is equal to the clique number of $\hat G_{p}$, thus $h_{p,2}=\alpha(G_p)$.
\end{proof}

\begin{prop} The independence number for the Cayley graphs $G_{2\ell}=\Cay(\mathbb Z_2^{2\ell}, S_{2\ell})$ is
\begin{equation}\label{eq:alpha}
\alpha(G_{2\ell})=2^{2\ell-1}-\frac{1}{2}\binom{2\ell}{\ell}.
\end{equation}
\end{prop}
\begin{proof}
The Cvetkovi\'c bound \cite{Cv} (see also \cite{CF}) for the independence number of a graph $G$ is $\alpha (G)\le p_0+\min\{p_+,p_-\}$, where $p_0$ is the multiplicity of $0$ as the eigenvalue of $G$, and $p_+$ and $p_-$ are the numbers of the positive and respectively negative eigenvalues of $G$.
From the formula for the eigenvalues of the Cayley graphs \cite{GR}, we get
$
\lambda (a) = \sum_{s\in S_{2\ell}}(-1)^{(a\cdot s)}=2\ell -2|a|+(-1)^{|a|},
$
for all $a\in \mathbb Z_2^{2\ell}$, where $|a|$ is the number of $1$ in $a$. Thus, all the eigenvalues are odd, thus $p_0=0$.  Hence $\lambda (a)>0$ if $0\le |a|\le \ell-1$, and  $\lambda (a)<0$ if $\ell+1\leq|a|\leq2\ell$. For $|a|=\ell$, $\lambda(a)$ is positive for even $\ell$ and it is negative for odd $\ell$. In both cases we get
$$
\min\{p_+,p_-\}=\sum_{k=0}^{\ell-1}\binom{2\ell}{k}.
$$

Thus,
$$
\alpha(G_{2\ell})\le\min\{ p_+,p_-\}=\sum_{k=0}^{\ell-1}\binom{2\ell}{k} =2^{2\ell-1}-\frac{1}{2}\binom{2\ell}{\ell}.
$$
It follows from  that $\alpha(G_{2\ell})\ge 2^{2\ell-1}-\frac{1}{2}\binom{2\ell}{\ell}$, thus, the formula \eqref{eq:alpha} is verified.
\end{proof}

\begin{prop} The independence number for the Cayley graphs $G_{2\ell+1}=\Cay(\mathbb Z_2^{2\ell+1}, S_{2\ell+1})$ is
\begin{equation}\label{eq:alpha1}
\alpha(G_{2\ell+1})=2^{2\ell}-\binom{2\ell}{\ell}.
\end{equation}
\end{prop}
\begin{proof}
We construct $G_{2\ell+1}=\Cay(\mathbb Z_2^{2\ell+1}, S_{2\ell+1})$  as two copies of $G_{2\ell}=\Cay(\mathbb Z_2^{2\ell}, S_{2\ell})$ placed above each other, with vertical edges connecting corresponding pairs of vertices from each layer. We want to prove that
$
\alpha(G_{2\ell+1})=2\alpha(G_{2\ell}).
$
Indeed, if $\alpha(G_{2\ell+1})>2\alpha(G_{2\ell})$, then at least one of the two congruent layers $L_1\cong L_2\cong G_{2\ell}$, say $L_1$ would contain at least  $\alpha(G_{2\ell})+1$ independent vertices of $G_{2\ell+1}$. Then, they would also form $\alpha(G_{2\ell})+1$ independent vertices of $G_{2\ell}$. That would be imply
$
\alpha(L_1)>\alpha(G_{2\ell}),
$
leading to the contradiction.

If $\alpha(G_{2\ell+1})<2\alpha(G_{2\ell})$, then by taking two disjoint maximal independent sets of  $G_{2\ell}$, and place one of them at the layer $L_1$ and another at $L_2$. We get $\alpha(G_{2\ell+1})\ge 2\alpha(G_{2\ell})$, leading again to the contradiction.

Thus, $\alpha(G_{2\ell+1})=2\alpha(G_{2\ell})$, verifying \eqref{eq:alpha1}, taking into account \eqref{eq:alpha}.
\end{proof}

We proved the following

\begin{thm}\label{th:h22} \begin{itemize}
\item[(i)] For $p=2\ell$, $h_{2\ell, 2}= 2^{2\ell-1}-\binom{2\ell}{\ell}/2$.
\item[(ii)] For $p=2\ell+1$, $h_{2\ell+1, 2}= 2^{2\ell}-\binom{2\ell}{\ell}$.
\end{itemize}
\end{thm}

\section{Total Complexity of Lasso with Normalized Data}\label{sec:comptotal}

The following theorem provides an upper bound on the total number of segments in a lasso with normalized data.

\begin{remark}
Let us recall a standard notation from the theory of hypergeometric functions (see e.g. \cite{AAR, VK}):
$$
{ }_{m}F_n(a_1,\dots, a_m;b_1\dots, b_n;z)=\sum_{j=0}^\infty\frac{(a_1)_j\dots(a_m)_j}{(b_1)_j\dots (b_n)_j}\frac{z^j}{j!},
$$
where $(a)_j=a(a+1)...(a+j-1)$. In particular, the Gaussian hypergeometric function is ${ }_{2}F_1(a,b;c;z)$
\end{remark}

\begin{thm}\label{thm:bound1} An upper bound on the total number of segments in a lasso with normalized data in dimension $p$
is
\begin{equation}\label{eq:ub}
\bar N(p)=\binom{p}{0}h_{p,2}+\binom{p}{1}h_{p-1,2}+...+\binom{p}{p-2}h_{2,2}+2,
\end{equation}
where $
h_{k,2}$ are given in Theorem \ref{th:h22}. Here,

\begin{equation}\label{eq:ubeven}
\bar N(2\ell)
=
\frac{
3^{2\ell}+4
-
{}_2F_1\left(
-\ell,\frac12-\ell;1;4
\right)
}{2}
-
2\ell\binom{2\ell-2}{\ell-1}
{}_3F_2\left(
\begin{matrix}
1-\ell,\;1-\ell,\;\frac12-\ell\\[1mm]
\frac32,\;\frac32-\ell
\end{matrix}
;\frac14
\right),
\end{equation}
and
\begin{equation}\label{eq:ubodd}
\bar N(2\ell+1)
=
\frac{
3^{2\ell+1}+4
-
{}_2F_1\left(
-\ell,-\ell-\frac12;1;4
\right)
}{2}
-
(2\ell+1)
{}_3F_2\left(
\begin{matrix}
-\ell,\;\frac12-\ell,\;\frac12\\[1mm]
1,\;\frac32
\end{matrix}
;4
\right).
\end{equation}
\end{thm}

In the above formulas, the infinite series that define hypergeometric functions ${}_2F_1$ and ${}_3F_2$ truncate into finite sums for the indicated sets of parameters (see also e.g. \cite{Dr}).

\begin{proof} There are $\binom{p}{p-k}$ ways to choose a $k$-dimensional face of a cross-polytope. According to Proposition \ref{prop:upperbp}, part (2), there are at most $h_{k,2}$ segments of a Lasso with normalized data at a given face of dimension $k$, for $k\ge 2$. However, according to Proposition \ref{prop:edgenorm}, there is at most one coordinate axis that contains a segment of a given solution of lasso with normalized data, which is always followed by the origin as the end of the solution.

To calculate $\bar N(p)$ in a closed form, we use the formula for $h_{p, 2}$ from Theorem \ref{th:h22}.
We show details of the computation for $p$ even, denoted as $p=2\ell$.
The formula \eqref{eq:ub} becomes
\begin{equation}\label{eq:suma}
\begin{aligned}
\bar N(2\ell)&=\binom{2\ell}{0}2^{2\ell-1}+\binom{2\ell}{1}2^{2\ell-2}+...+\binom{2\ell}{2\ell-2}2^1+2\\
&-\frac{1}{2}\Bigg(\binom{2\ell}{0}\binom{2\ell}{\ell}+\binom{2\ell}{2}\binom{2\ell-2}{\ell-1}+\binom{2\ell}{4}\binom{2\ell-4}{\ell-2}\dots\Bigg)\\
&-\Bigg(\binom{2\ell}{1}\binom{2\ell-2}{\ell-1}+\binom{2\ell}{3}\binom{2\ell-4}{\ell-2}+\binom{2\ell}{5}\binom{2\ell-6}{\ell-3}\dots\Bigg)\\
&=N_1(2\ell)-\frac12N_2(2\ell)-N_3(2\ell).
\end{aligned}
\end{equation}
We consider the three sums in \eqref{eq:suma} separately, starting from the first one:
$$
\begin{aligned}
N_1(2\ell)&=\frac{1}{2}\Bigg(\binom{2\ell}{0}2^{2\ell}+\binom{2\ell}{1}2^{2\ell-1}+...+\binom{2\ell}{2\ell-1}2+\binom{2\ell}{2\ell}\Bigg)+2-2\ell-\frac12\\
&=\frac{3^{2\ell}+3}{2}-2\ell.
\end{aligned}
$$
To calculate  $N_2(2\ell)$ we use the identity for nonegative integers $\ell$ and $m$, $\ell\geq m$:
\begin{equation}\label{eq:identity}
\binom{2\ell}{2m}\binom{2\ell-2m}{\ell-m}=\binom{2\ell}{\ell-m}\binom{\ell+m}{\ell-m}.
\end{equation}
Using \eqref{eq:identity}, the sum $N_2(2\ell)$ is rewritten as
\begin{equation}\label{eq:N2}
N_2(2\ell)=\binom{2\ell}{0}\binom{2\ell}{0}+\binom{2\ell}{1}\binom{2\ell-1}{1}+...+\binom{2\ell}{\ell}\binom{\ell}{\ell}=\sum_{k=0}^\ell \binom{2\ell}{k}\binom{2\ell-k}{k}.
\end{equation}
\begin{remark} The expression \eqref{eq:N2} gives a nice combinatorial interpretation of $N_2$. Let us consider $1\times 2\ell$ strip formed by $2\ell$ unit squares in a line. The sum $N_2$ is the number of different ways to color the strip in three colors red, blue and green, with the condition that there is the same number of red and blue squares.
\end{remark}
Using the notation of George Andrews (see e.g. \cite{An, AB}) , the trinomial coefficient $\binom{n}{k}_2$ is the coefficient of $x^{n+k}$ in the expansion of $(1+x+x^2)^n$, where $n\geq0$ and $-n\leq k\leq n$ and we get:
$$N_2(2\ell)=\binom{2\ell}{0}_2={}_2F_1(-\ell,\frac{1}{2}-\ell;1;4).$$
A similar calculation gives $N_3(2\ell)$.

The formula for $\bar N(p)$ for odd $p$ is derived in an analogous manner.
\end{proof}

The first values of $\bar N$ are
$
\bar N(1)=2,
\bar N(2)=3,
\bar N(3)=7,
\bar N(4)=21, \bar N(5)=67,\bar N(6)=214, \dots.$

\begin{prop}\label{prop:totalp3}
The sharp upper bound for the total complexity of a lasso with normalized data for $p=3$ is $6$, that is smaller than $\bar N(3)=7$.
\end{prop}
\begin{proof}
 Since $h_{3,2}=2$ (according to Theorem \ref{th:h22}), we see that in a $p=3$ normalized lasso, one can have only two segments from $3$-dimensional orthants, they can be separated by only one edge from a $2$-dimensional orthant, leading to three segments in total. Once the lasso reached the second segment of  a $3$-dimensional orthant, it can make not more than two segments before getting to the origin. Thus, the sharp upper bound for the total complexity a lasso with normalized data for $p=3$ is at most $6$.

Example \ref{exam:transfer3d} provides a case that shows that the total number of $6$ segments can be achieved by a lasso with normalized data for $p=3$. See also Figure \ref{sl:4.5R}. The sixth segment is the origin of the $3$-dimensional space, as the degenerate segment that ends the lasso.

Thus, the sharp upper bound for the total complexity a lasso with normalized data for $p=3$ is exactly $6$.

\end{proof}

\begin{dfn}\label{def:code} We call  \emph{a normalized lasso-type code of dimension $p$} any finite sequence  of $p$-tuples of symbols $\{-1, 0, 1\}$, that satisfies the following axioms: (i) every two consecutive $p$-tuples differ by exactly one change, where the changes are $1\leftrightarrow 0$ or $-1\leftrightarrow 0$; (ii) one $p$-tuple can appear at most once in the code; (iii) for any given $k<p$, any two $p$-tuples with exactly $k$ zeros, that occupy the same spots, have at least two remaining nonzero spots that match; (iv) there is only one $p$-tuple in the code that has all but one spots equal to zero; that $p$-tuple is followed by the zero $p$-tuple, which is the last $p$-tuple in the code. We refer to the $p$-tuples in a code also as $p$-spot labels.
\end{dfn}

\begin{thm}\label{thm:bound3} An upper bound for the number of labels in a normalized lasso-type code of dimension $p$ is $\hat N(p)$, where
\begin{align}
\hat N(2\ell+1)&=2\Big(\sum_{j=0}^{\ell-1}\binom{2\ell+1}{2j}h_{2(\ell-j)+1,2}+1\Big)+1,\\
\hat N(2\ell)&=2\Big(\sum_{j=1}^{\ell-1}\binom{2\ell}{2j-1}h_{2\ell-2j+1,2}+1\Big)+1.
\end{align}
\end{thm}
\begin{proof} The proof follows from Theorem \ref{thm:bound1}. We apply a similar consideration  as in  relation \eqref{eq:ub}, taking into account that after a label with odd number of zeros it goes a label with even number of zeros and vice versa, according to axiom (i) from Definition \ref{def:code}.
\end{proof}

\begin{prop} The closed formulas for $\hat N(p)$ are:
\begin{align}
\hat N(2\ell+1)&=\frac{3^{2\ell+1}+7}{2}
-
2(2\ell+1)\,
{}_3F_2\!\left(
\begin{matrix}
\frac{1}{2},\;-\ell,\;\frac{1}{2}-\ell\\[1mm]
\frac{3}{2},\;1
\end{matrix}
;8
\right),\\
\hat N(2\ell)&=\frac{3^{2\ell}+5}{2}
-
4\ell\,
{}_3F_2\!\left(
\begin{matrix}
\frac{1}{2},\;1-\ell,\;\frac{1}{2}-\ell\\[1mm]
\frac{3}{2},\;1
\end{matrix}
;4
\right).
\end{align}

\end{prop}
\begin{proof} The proof follows the lines of the proof of Theorem \ref{thm:bound1}, using, for example in the $p=2\ell$ cases:
\begin{align*}
\sum_{j=1}^{\ell-1}\binom{2\ell}{2j-1}2^{2(\ell-j)}&=\frac{3^{2\ell-1}}{4}-2\ell,\\
\sum_{j=1}^{\ell-1}\binom{2\ell}{2j-1}\binom{2(\ell-j)}{\ell-j}&=2\ell\Big({}_3F_2\left(
\begin{matrix}
\frac12,\;\frac12-\ell,\;1-\ell\\[1mm]
1,\;\frac32
\end{matrix}
;4
\right)-1\Big),
\end{align*}
with similar formulas in $p=2\ell+1$ cases.
\end{proof}

\begin{exm} A normalized lasso-type code of dimension $p=4$ with exactly $\hat N(4)=19$ lines:
$$
\begin{aligned}
&(1,1,1,1), (0,1,1,1), (-1,1,1,1), (-1,0,1,1), (0,0,1,1),(1,0,1,1), (1, -1, 1,1),\\
&(1,-1,0,1), (1,0,0,1), (1,1,0,1), (1,1,-1,1), (1,1,-1,0),(1,1,0,0), (1, 1, 1,0),\\
&(1,1,1,-1), (0,1, 1,-1), (0,1,0,-1), (0,0,0,-1), (0, 0, 0,0).
\end{aligned}
$$
\end{exm}

From the previous examples we see that  $\hat N(4)=19$ is the sharp upper bound for the number of $p$-spot labels in a normalized lasso-type code of dimension $p=4$. For $p=3$, such a sharp upper bound is $6$.

We conclude the paper, by formulating  the following questions:

\begin{itemize}
 \item[(Q1)] What is $h_{p,k}$ for all $p$ and all $p-2>k>2$? ($h_{p, p-1}=2$, $h_{p, p-2}=p+1$.)
 \item[(Q2)] What is the sharp upper bound for the number of $p$-spot labels in a normalized lasso-type code of dimension $p$ for $p\ge 5$? Is this $\hat N(p)$? (For $p=5$, we have an example showing that it indeed is $\hat N(5)=63$.)
  \item[(Q3)] What is the  sharp upper bound for the total  complexity of a lasso with normalized data in dimension $p$, for $p>4$? Is this $\hat N(p)$?
  \item[(Q4)] Does for every normalized lasso-type code there exist a  lasso with normalized data that generates it?
  \end{itemize}

\subsection*{Acknowledgements}
We are deeply grateful  to Pankaj Choudhary for his constant support, encouragement, and interesting discussions.  We are also grateful to Robert Tibshirani, Herbert Edelsbrunner, B\'ela Bollob\'as, and Nathan Williams for useful discussions.
This research has been partially supported by the Simons Foundation grant no. 854861.


\begin{thebibliography}{99}

\bibitem {An} Andrews, G. (1990) Euler's 'exemplum memorabile inductionis fallacis' and q-Trinomial Coefficients, \emph{J. Amer. Math. Soc.}, \textbf{3}, 653-669.

\bibitem {AB}  Andrews, G. E., Baxter R. J., (1987) Lattice gas generalization o/the hard hexagon model. III.
q-trinomial coefficients, \emph{J. Statist. Phys.}, \textbf{47}, 297-330.

\bibitem{AAR} Andrews, G., Askey, R., Roy, R., (1999) \emph{Special Functions}, Encyclopedia of Mathematics and ist Applications 71, Cambridge University Press


 \bibitem{CF} Calderbank, A., R. Frankl P., (1992) Improved upper bounds concerning the
Erd\"os-Ko-Rado theorem, \emph{Combin. Probab. Comput.} \textbf{1}, no. 2, 115--122.



\bibitem{Co} Coxeter, H. S. M. (1973) \emph{Regular Polytopes}, Dover.

\bibitem{Cr}Cram\'er, H. (1946) \emph{Mathematical Methods of Statistics,} Princeton: Princeton University Press.

\bibitem{Cv} Cvetkovi\' c D. M.  (1971) Graphs and their spectra, \emph{Univ. Beograd. Publ.
Elektrotehn. Fak. Ser. Mat. Fiz.} no. 354, 1--50.

\bibitem{Dr} Dragovi\'c, V. (2002) The Appell hypergeometric functions and separable mechanical systems, \emph{Journal of Physics A: Mathematical and General},  Vol. 35, No. 9, p. 2213--2222.

\bibitem{DG2022a} Dragovi\'c, V. Gaji\'c, B. (2023) Points with rotational ellipsoids of inertia, envelopes of hyperplanes which equally fit the system of points in $R^k$, and ellipsoidal billiards, \emph{Physica D: Nonlinear Phenomena}, 15 p. Volume 451, 133776


\bibitem{DG2023} Dragovi\'c, V. Gaji\'c, B. (2025) Orthogonal and Linear Regressions and Pencils of Confocal Quadrics, \emph{Statisical Science}, Vol. 40, No. 2, 289 -- 312.

\bibitem{DG2023a} Dragovi\'c, V., Gaji\'c, B. (2025). Supplement to ``Orthogonal and Linear Regressions and Pencils of Confocal Quadrics'', DOI:10.1214/24-STS938SUPP.


\bibitem{Ge} Gentle, J. E. (1998). \emph{ Numerical Linear Algebra for Applications in Statistics}. Springer.




\bibitem{GR}
Godsil , C. and Royle G. (2001) \emph{Algebraic Graph Theory}, Springer.




\bibitem{GSST} Graczyk, P., Schneider, U., Skalski, T. and Tardivel, P. (2023) A Unified Framework for Pattern Recovery in Penalized and Thresholded Estimation and its Geometry; arxiv:2307.10158[math.ST].


\bibitem{HJ} Horn, R. A., Johnson, C. R. (1985). \emph{Matrix Analysis.} Cambridge University Press.





 \bibitem{FMF}
   Friendly, M., Monette, G., Fox, J. (2013). Elliptic insights: understanding statistical methods through elliptical geometry. \emph{Statist. Sci.} \textbf{28}, no. 1, 1--39.




\bibitem{HTF} Hastie, T.,   Tibshirani, R.,  Friedman, J. (2009) \emph{The elements of statistical learning} Second Edition, Springer

\bibitem{HTW} Hastie, T.,   Tibshirani, R.,  Wainwright, M. (2015) \emph{Statistical Learning with Sparsity:
The Lasso and Generalizations}, CRC Press


\bibitem{KM} Klee, V., Minty, G. J., (1972) How good is the simplex algorithm? In Shisha, O (ed.) \emph{Inequalities, III} pp. 159-175, Academic Press, New York.

\bibitem{Kl} Kleitman, D., (1966) On a combinatorial conjecture of Erd\"os, \emph{Journal of Combinatorial Theory}, Vol. 1, Issue 2, p. 209--214.

\bibitem{MY} Mairal J.,  Yu B. (2012) Complexity Analysis of the Lasso Regularization Path, \emph{Proceedings of the International Conference on Machine Learning 29}.


\bibitem{MF} Matejka, J.,  Fitzmaurice, G. (2017) Same Stats, Different Graphs:Generating Datasets with Varied Appearance andIdentical Statistics through Simulated Annealing, \emph{CHI '17: Proceedings of the 2017 CHI Conference on Human Factors in Computing SystemsMay}, 1290–1294,  DOI: http://dx.doi.org/10.1145/3025453.3025912

\bibitem{Ro} Rogers, A. J. (2013). Concentration Ellipsoids, Their Planes of
Support, and the Linear Regression Model, \emph{Econometric Reviews}, 32:2, 220-243, DOI:10.1080/07474938.2011.608055.


\bibitem {ST} Schneider, U., Tardivel, P (2022) The Geometry of Uniqueness, Sparsity and Clustering in
Penalized Estimation, \emph{Journal of Machine Learning Research} \textbf{23}, 1--36

\bibitem{SKM} Stamey, T., Kabalin, J., McNeal, J., Johnstone, I., Freiha, F., Redwine, E. and Yang, N. (1989) Prostate specific antigen in the diagnosis and treatment of adenocarcinoma of the prostate II radical prostatectomy treated patients, \emph{Journal of Urology} 16: 1076--1083.


\bibitem{Tib} Tibshirani, Robert. (1996), Regression shrinkage and selection via the lasso, \emph{Journal of the Royal Statistical Society: Series B} 58(1), 267--288
\bibitem{Tib2} Tibshirani, Ryan. (2013), The Lasso problem and uniqueness, \emph{Electronic Journal
of Statistics}, 7, 1456--1490.


\bibitem{VK} Vilenkin, N. and Klimyk A, (1995) \emph{Representation of Lie groups and special functions. Recent Advances},
Dordrecht: Kluwer, p 497.

\bibitem{Zi} Ziegler, G. M. (1995) \emph{Lectures on Polytopes}, GTM, V. 152, Springer.

\end{thebibliography}
\end{document}